\documentclass{article}
\usepackage{amssymb, amsmath}
\usepackage{amsthm}
\newtheorem{lem}{\bf Lemma}[section]%
\newtheorem{thm}[lem]{\bf Theorem}%
\newtheorem{prop}[lem]{\bf Proposition}%
\newtheorem{cor}[lem]{\bf Corollary}%
\newtheorem{rem}[lem]{\sc Remark}%
{}%
\newenvironment{pf}{\noindent {\sf Proof.}}{$\Box$\par\medskip}%
\newcommand{\be}{\begin{equation}}
\newcommand{\ee}{\end{equation}}
\newcommand{\ol}{\overline}

\newcommand{\wh}{\widehat}

\begin{document}


 \title{Parity, eulerian subgraphs and the Tutte~polynomial}
 \author{A.J. Goodall}
\maketitle

\begin{abstract}
Identities obtained by elementary finite Fourier analysis are used to
derive a variety of evaluations of the Tutte polynomial of a graph $G$
at certain points $(a,b)$ where $(a-1)(b-1)\in\{2,4\}$. 
These evaluations are expressed in terms of eulerian subgraphs of
$G$ and the size of subgraphs modulo $2,3,4$ or $6$. In particular, a
graph is found to have a nowhere-zero
$4$-flow if and only if there is a correlation between the event that
three subgraphs $A,B,C$ chosen uniformly at random have
pairwise eulerian symmetric differences and the event that
$\lfloor\frac{|A|+|B|+|C|}{3}\rfloor$ is even. 
Some further evaluations of the Tutte polynomial at points $(a,b)$
where $(a-1)(b-1)=3$ are also given that illustrate the unifying power
of the methods used. The connection between results of Matiyasevich
\cite{M04}, Alon and Tarsi \cite{AT92} and Onn \cite{Onn04}
is highlighted by indicating how they may all be derived by the techniques adopted in
this paper. \end{abstract}




\section{Introduction}

Let $G$ be a finite graph, loops and parallel edges permitted. This article continues a series of papers
\cite{AJG05}, \cite{CCC07} using elementary finite
Fourier analysis to derive evaluations of the Tutte polynomial
$T(G;x,y)$ at certain points $(x,y)=(a,b)$ where
$(a-1)(b-1)\in\{2,3,4\}$. Letting $H_q$ denote the hyperbola
$\{(x,y):(x-1)(y-1)=q\}$,
the evaluations of the Tutte polynomial on $H_2$ and $H_4$ are expressed in terms of eulerian subgraphs of
$G$ and the size of subgraphs modulo $2,3,4$ or $6$, the evaluations
on $H_3$ in terms of directed
eulerian subgraphs and size modulo $3$. 

The theorems
obtained in this article are analogous to the well-known fact that the
eulerian subgraphs of $G$ all have the same size modulo $2$
if $G$ is bipartite and otherwise the number
of eulerian subgraphs of odd size is equal to the number of even size. They seem however to be more
elusive of explanation than this straightforward example. 
In particular, Theorem~\ref{sum cubes with q=6}  states that a
graph $G$ has a nowhere-zero
$4$-flow if and only if there is a correlation between (i) the event that
three subgraphs $A,B,C$ chosen uniformly at random have
pairwise eulerian symmetric differences and (ii) the event that
$\lfloor\frac{|A|+|B|+|C|}{3}\rfloor$ is even. A companion to this
theorem is Theorem~\ref{squares with q=4} that $G$ is eulerian (has a
nowhere-zero $2$-flow) if and only if there is a correlation between (i)
the event that subgraphs $A,B$ chosen uniformly at random have
eulerian symmetric difference and (ii) the event that
$\lfloor\frac{|A|+|B|}{2}\rfloor$ is even. 
In this language, a graph $G$ is bipartite (has a nowhere-zero $2$-tension) if
and only if there is a correlation between
the event that a subgraph $A$ chosen uniformly at random is
eulerian and the event that
$|A|$ is even.

In section~\ref{poset} we develop a result of Onn \cite{Onn04}
which states a parity criterion for the existence of nowhere-zero
$q$-flows of a graph.
 This leads to another correlation between a parity event and an
event involving eulerian subgraphs of $G$, this time connected to
whether or not $G$ has a proper vertex $4$-colouring. 
Onn in his paper uses
the algebraic method used by Alon and Tarsi \cite{AT92}, \cite{AT97}, \cite{A99} in their proof of a parity criterion for the existence of
proper $q$-colourings of $G$ in terms of eulerian subdigraphs of an
orientation of $G$.  

In section~\ref{cubic graphs} we consider graphs with a cycle
double cover by triangles (such as plane triangulations) and derive a
criterion for the existence of a proper vertex $4$-colouring of $G$ in
terms of the correlation between (i) the event that a pair of subgraphs
$A,B$ of $G$ are eulerian and between them cover all the edges of $G$
and (ii) the event that $|A|\equiv|B|\pmod 3$. 

In section~\ref{4 regular} we prove theorems of a similar
character to the cited results of Alon and Tarsi, involving evaluations of
the Tutte polynomial on $H_3$, this time confining our attention to
$4$-regular graphs (such as line graphs of cubic graphs). These results stem also from the work of Matiyasevich \cite{M04}, whose probabilistic restatements of the Four
Colour Theorem inspired the mode of expression for the Tutte
polynomial evaluations throughout this paper. 

The method of proof throughout is to use identities from
elementary Fourier analysis from which the
interpretations of the Tutte polynomial evaluations can be extracted. 
The results of section~\ref{H2 H4} (the main theorems of which are to
be found in sections~\ref{square root of unity}
to \ref{sixth root of unity}) ultimately derive from  Lemma~\ref{sum
  cubes}, which comprises a set of identities which can be found in \cite{E98} and
\cite{Mi05}. In the language of coding theory, these identities relate the sum of powers of coset weight enumerators of a
binary code to the Hamming weight enumerator of the code. In the
context of this paper, the code is the cycle space of a graph. In
sections~\ref{poset}, \ref{cubic graphs} and \ref{4 regular} the main
tool is the discrete version of the Poisson summation formula (or, in the context of
coding theory, the MacWilliams duality
theorem for complete weight enumerators). A
simultaneous generalisation of Lemma~\ref{sum cubes} and the Poisson
summation formula is presented in Lemma~\ref{cwe sum}: this is used in
Section~\ref{4 regular}.

 A more expansive
exposition of the material in
sections~\ref{sec:fourier} and \ref{hwes} of the present article can be found in
\cite{CCC07} but is included here for convenience. All the facts quoted
without reference in section~\ref{sec:fourier} on Fourier transforms can be found for example in
\cite{Terras99}.

Eulerian subgraphs of $G$ are cycles in the graphic matroid
underlying $G$. In \cite{Kung07} the results of \cite{AJG05} are extended from graphs
to matrices. In a similar way, the results of the present paper depend only on
the cycle space of $G$ and could be extended to binary matroids. 



\subsection{Notation and definitions}

Let $G=(V,E)$ be a graph. In section~\ref{H2 H4} we consider $\mathcal{C}_2\subseteq\mathbb{F}_2^E$ the subspace of
  $2$-flows (eulerian subgraphs, cycles) of $G$, and
  $\mathbb{F}_2^E/\mathcal{C}_2$ the set of cosets of $\mathcal{C}_2$ in the additive group
  $\mathbb{F}_2^E$. 
The quotient space $\mathbb{F}_2^E/\mathcal{C}_2$ is isomorphic to the orthogonal subspace
$\mathcal{C}_2^\perp$ of $2$-tensions (cutsets, cocycles) of $G$.
 
The rank of $G$ is defined by $r(G)=|V|-k(G)$, where $k(G)$ is the
number of components of $G$, and the nullity by $n(G)=|E|-r(G)$. The
subspace $\mathcal{C}_2$ has dimension $n(G)$ over $\mathbb{F}_2$, and
$\mathcal{C}_2^\perp$ dimension $r(G)$.

To each subset $A\subseteq E$ there is a subgraph $(V,A)$ of $G$
obtained by deleting the edges in $E\setminus A$ from $G$. For short
this subgraph will be referred to just by its edge set $A$,
 ``the
subgraph $A$'' meaning the graph $(V,A)$.  
A subgraph $A$ is {\em eulerian} if all its vertex degrees are even.
The subspace $\mathcal{C}_2$ of $2$-flows of $G$ may be identified with the
set of eulerian subgraphs of $G$.
For two subgraphs $A,B$ of $G$ the
symmetric difference $A\bigtriangleup B$ corresponds to 
addition of the indicator vectors of $A$ and $B$ in $\mathbb{F}_2^E$. The size $|A|$ of
$A$ is equal to the Hamming weight of the indicator vector of $A$.
Two subsets $A,B\subseteq E$ belong to the same coset of $\mathcal{C}_2$
if and only if $A\bigtriangleup B$ is eulerian, and this is the case
if and only if the subgraphs $A,B$ have the same degree sequence
modulo $2$.

 The space of $\mathbb{F}_4$-flows of $G$ will be denoted by
 $\mathcal{C}_4$. The space $\mathcal{C}_4$ is isomorphic to $\mathcal{C}_2\times\mathcal{C}_2$. From this observation a graph
  $G=(V,E)$ has a nowhere-zero $4$-flow if and only if there are two eulerian
  subgraphs which together cover $E$, whence the well-known
  equivalence of the Four Colour Theorem with existence of an edge
  covering of any given planar graph by two of its eulerian subgraphs.

 For $A\subseteq E$ the rank $r(A)$ of $A$
 is defined to be the rank of the subgraph $(V,A)$. The {\em Tutte polynomial} of
 $G$ is defined by
$$T(G;x,y) =\sum_{A\subseteq E}(x-1)^{r(E)-r(A)}(y-1)^{|A|-r(A)}.$$
The hyperbolae $H_q=\{(x,y)\,:\,(x-1)(y-1)=q\}$ for $q\in\mathbb{N}$ play a special role in the theory
of the Tutte polynomial, summarised for example in \cite[\S
3.7]{DW93}. 
In particular,  $(-1)^{r(G)}q^{k(G)}T(G;1-q,0)=P(G;q)$ is
 the number of proper vertex $q$-colourings of $G$ and
 $(-1)^{n(G)}T(G;0,1-q)=F(G;q)$ is the number of nowhere-zero $\mathbb{Z}_q$-flows of
 $G$. 

On $H_2$ there are the evaluations
 $P(G;2)=(-1)^{r(G)}2^{k(G)}T(G;-1,0)$ and $F(G;2)=(-1)^{n(G)}T(G;0,-1)$. The Tutte
 polynomial on $H_2$ is the partition function of the Ising model of
 statistical physics, or, what is the same thing, the Hamming weight
 enumerator of the subspace $\mathcal{C}_2$ of $2$-flows (this is Van
 der Waerden's eulerian expansion of the Ising model~\cite{vdW41}). In
 section~\ref{H2 H4} the Tutte polynomial at the points
 $(-2,\frac{1}{3})$ and $(-\frac{1}{2},-\frac{1}{3})$ on $H_2$, in addition to $(-1,0)$ and $(0,-1)$, are
 given an interpretation in terms of eulerian subgraphs of $G$.

On $H_4$ we have $P(G;4)=(-1)^{r(G)}4^{k(G)}T(G;-3,0)$ and 
 $F(G;4)=(-1)^{n(G)}T(G;0,-3)$. Also, $T(G;-1,-1)=(-1)^{|E|}(-2)^{\dim(\mathcal{C}_2\cap\mathcal{C}_2^\perp)}$, where
 $\mathcal{C}_2\cap\mathcal{C}_2^\perp$ is the bicycle space of $G$~\cite{RR78}. 
 The Tutte
 polynomial on $H_4$ is the partition function of the $4$-state Potts
 model, which coincides with the Hamming weight
 enumerator of the subspace $\mathcal{C}_2\times\mathcal{C}_2$ of
 $\mathbb{F}_4$-flows. Evaluations of the Tutte polynomial at the point
 $(-2,-\frac{1}{3})$ on $H_4$ as well as $(-3,0)$, $(0,-3)$ and
 $(-1,-1)$ are also given an interpretation in section~\ref{H2 H4} in terms of eulerian
 subgraphs of $G$.

By MacWilliams duality (see section~\ref{hwes} below), the
interpretions that we give for evaluations of the Tutte polynomial at points $(a,b)$ in terms of eulerian
subgraphs (cycles, $\mathcal{C}_2$) become interpretations for points $(b,a)$ in terms of
cutsets (cocycles, $\mathcal{C}_2^\perp$). In the same way, for
example, there is the expansion of the
 Ising model of a graph over its cutsets (bipartitions)~\cite[\S 4.3]{DW93}
corresponding to Van der Waerden's expansion over
eulerian subgraphs.

\subsection{The Fourier transform}\label{sec:fourier}

In this section we summarise the facts about the Fourier transform on rings such as
$\mathbb{F}_q$ and $\mathbb{Z}_q$ (the integers modulo $q$) that the reader needs to be aware of
in this paper. 

Let $Q$ be a commutative ring (either $\mathbb{F}_q$ or
$\mathbb{Z}_q$ in the sequel) and
$Q^E$ the set of all vectors $x=(x_e:e\in E)$ with entries in $Q$
indexed by $E$. The indicator function $1_k$ for $k\in
  Q$ is defined
  by $1_k(\ell)=1$ if $k=\ell$ and $0$ otherwise. A subset $S\subseteq
  Q$ has indicator function defined by $1_S=\sum_{k\in S}1_k$.

A {\em character} of $Q$ is a homomorphism $\chi:Q\rightarrow\mathbb{C}^\times$ from the additive
group $Q$ to the multiplicative group of $\mathbb{C}$. The set of
characters form a group $\wh{Q}$ under pointwise
multiplication isomorphic to $Q$ as an abelian group. For $Q^E$ the group of characters
$\wh{Q^E}$ is isomorphic to $\wh{Q}^E$. 

 For each $k\in Q$, write $\chi_k$ for the image of $k$ under a fixed
isomorphism of $Q$ with $\wh{Q}$. In particular, the principal
(trivial) character $\chi_0$ is defined by $\chi_0(\ell)=1$ for all $\ell\in
Q$, and $\chi_{-k}(\ell)=\ol{\chi_k(\ell)}$ for all $k,\ell\in Q$, where the bar
denotes complex conjugation. 
A character $\chi\in\wh{Q}$ is a {\em generating character} for $Q$ if
$\chi_k(\ell)=\chi(k\ell)$ for
each character $\chi_k\in\wh{Q}$. The ring $\mathbb{Z}_q$ 
has a generating character $\chi$ defined by $\chi(k)=e^{2\pi i
  k/q}$. (The fixed
isomorphism $k\mapsto \chi_k$ of $\mathbb{Z}_q$ with $\wh{\mathbb{Z}_q}$ in this case is
given by taking $\chi_k(\ell)= e^{2\pi ik\ell/q}$.) 
The field $\mathbb{F}_q$ with $q=p^m$ has a generating
character $\chi$ defined by $\chi(k)=e^{2\pi i
\mbox{\rm \tiny Tr}(k)/p}$, where ${\rm Tr}(k)=k+k^p+\cdots +
k^{p^{m-1}}$ is the trace of $k$. (Here the isomorphism
$\mathbb{F}_q\rightarrow\wh{\mathbb{F}_q},\; k\mapsto\chi_k$ is given
by taking $\chi_k(\ell)=e^{2\pi i\mbox{\rm \tiny Tr}(k\ell)/p}$.)
If $\chi$ is a
generating character for $Q$ then $\chi^{\otimes E}$, defined by
$\chi^{\otimes{E}}(x)=\prod_{e\in E}\chi(x_e)$ for $x=(x_e:e\in E)\in Q^E$, is a generating character for $Q^E$. The euclidean
inner product (dot product) on $Q^E$ is defined by $x\cdot y=\sum x_ey_e$. Since
$\prod\chi(x_e)=\chi(\sum x_e)$, it follows that
$\chi^{\otimes E}$ satisfies $\chi^{\otimes E}_x(y)=\chi(x\cdot y)$ for $x,y\in
Q^E$. 

The vector space
$\mathbb{C}^{Q^E}$ over $\mathbb{C}$ of all functions from $Q^E$ to
$\mathbb{C}$ is an inner product space with Hermitian
inner product $\langle\, ,\,\rangle$ defined for
$f,g:Q^E\rightarrow\mathbb{C}$ by
$$\langle f,g\rangle=\sum_{x\in Q^E}f(x)\ol{g(x)}.$$

Fix an isomorphism $k\mapsto \chi_k$ of $Q$ with $\wh{Q}$ and let
$\chi$ be a generating character for $Q$ such that
$\chi_k(\ell)=\chi(k\ell)$.
For $f\in\mathbb{C}^{Q}$ the {\em Fourier transform}
$\wh{f}\in\mathbb{C}^{Q}$ is defined for $k\in Q$ by 
\[\wh{f}(k)=\langle f,\chi_k\rangle=\sum_{\ell\in
  Q}f(\ell)\chi(-k\ell).\]

The Fourier transform of a function $g:Q^E\rightarrow\mathbb{C}$ is
then given by
$$\wh{g}(x)=\sum_{y\in Q^E}g(y)\chi(-x\cdot y).$$

In the space $\mathbb{C}^{Q^E}$, the Fourier inversion
  formula is $\wh{\wh{f}}(x)=q^{|E|}f(-x)$, or
$$f(y)=q^{-|E|}\langle \wh{f},\chi_{-y}\rangle=q^{-|E|}\sum_{x\in
    Q^E}\wh{f}(x)\chi(x\cdot y).$$
 Plancherel's or Parseval's 
  identity is $\langle f,g\rangle=q^{-|E|}\langle\wh{f},\wh{g}\rangle.$

For a subset $\mathcal{S}$ of $Q^E$, the annihilator
$\mathcal{S}^{\sharp}$ of
$\mathcal{S}$ is defined by
$\mathcal{S}^{\sharp}=\{y\in Q^E:\forall_{x\in
  \mathcal{S}}\,\,\chi_x(y)=1\}.$ 
 The annihilator $\mathcal{S}^\sharp$ is a subgroup of $Q^E$
 isomorphic to $Q^E/\mathcal{S}$. When $\mathcal{S}$ is a
 $Q$-submodule of $Q^E$ and $Q$ has a generating character,
 the annihilator of $\mathcal{S}$ is equal to the orthogonal $\mathcal{S}^\perp$
to $\mathcal{S}$ (with respect
to the euclidean inner product), defined by
$\mathcal{S}^\perp=\{y\in Q^E:\forall_{x\in \mathcal{S}}\,\, x\cdot y=0\}.$

A key property of the Fourier transform is that for a $Q$-submodule
$\mathcal{S}$ of $Q^E$
$$\wh{1_\mathcal{S}}(y)=\sum_{x\in \mathcal{S}}\chi_x(y)=|\mathcal{S}|1_{\mathcal{S}^\sharp}(y),$$
and the Poisson summation formula is that
\[\sum_{x\in \mathcal{S}}f(x+z)=\frac{1}{|\mathcal{S}^\sharp|}\sum_{x\in
    \mathcal{S}^\sharp}\wh{f}(x)\chi_{z}(x).\]

For $Q$-submodule $\mathcal{S}$ of $Q^E$ the coset
$\{x+z:x\in\mathcal{S}\}$ of $\mathcal{S}$ in the additive group $Q^E$
is denoted by $\mathcal{S}+z$, an element
of the quotient module $Q^E/\mathcal{S}$.

\subsection{Flows, tensions and Hamming weight enumerators}\label{hwes}

For the moment we continue with $Q$ a commutative ring on $q$ elements with a generating
character. Let $\mathcal{S}$ be a subset of $Q^E$, the set
of vectors with entries in $Q$ indexed by edges of $G$. 
The Hamming weight of a vector $x=(x_e:e\in E)\in Q^E$ is defined by
$|x|=\#\{e\in E:x_e\neq 0\}$. The Hamming weight enumerator of
$\mathcal{S}$ is defined by
$${\rm hwe}(\mathcal{S};t)=\sum_{x\in\mathcal{S}}t^{|E|-|x|},$$
the exponent being the number of zero entries in $x$.

When $\mathcal{S}$ is a $Q$-submodule of $Q^E$ the MacWilliams
  duality theorem states that 
\be\label{MacWilliams}{\rm
  hwe}(\mathcal{S};t)=\frac{(t-1)^{|E|}}{|\mathcal{S}^{\perp}|}{\rm
  hwe}\left(\mathcal{S}^\perp;\frac{t+q-1}{t-1}\right),\ee
which follows from the Poisson summation formula with
  $f=(t1_0+1_{Q\setminus 0})^{\otimes E}$. 

A {\em $Q$-flow} of $G$ is defined with reference to a ground
orientation $\gamma$ of $G$; the number of $Q$-flows of a given Hamming
weight is independent of $\gamma$. A vector $x\in Q^E$ is a $Q$-flow of $G$ if, for each
vertex $v\in V$, 
$$\sum_{e\in E}\gamma_{v,e}x_e=0,$$
where $\gamma_{v,e}=+1$ if $e$ is directed into $v$, $\gamma_{v,e}=-1$
is $e$ is directed out of $v$, and $\gamma_{v,e}=0$ if $e$ is not
incident with $v$. The $Q$-flows form a $Q$-submodule of $Q^E$ whose
orthogonal $Q$-submodule is the set of $Q$-tensions of
$G$. The latter comprise the set
of $y\in Q^E$ such that there exists a
vertex $Q$-colouring $z\in Q^V$ with $y_e=z_u-z_v$ for
all edges $e=(u,v)$ (directed by the orientation $\gamma$). To each
$Q$-tension $y$ there correspond $q^{k(G)}$ vertex $Q$-colourings for
which $y_e=z_u-z_v$. 

A {\em nowhere-zero} $Q$-flow has Hamming weight $|E|$, and likewise a
nowhere-zero $Q$-tension. To a nowhere-zero $Q$-tension $y$
corresponds a set of $q^{k(G)}$ proper vertex $Q$-colourings of $G$, i.e., $z\in Q^V$ such that $z_u\neq z_v$ whenever $u$ is adjacent to $v$ in $G$.
 
The Tutte polynomial on the hyperbola $H_q$ is related to the Hamming weight
enumerator of the set of $Q$-flows via the identity
\be\label{qflows tutte}
{\rm hwe}(\mbox{\rm $Q$-flows of
  $G$};t)=(t-1)^{n(G)}T\left(G;t,\frac{t+q-1}{t-1}\right).
\ee
By \eqref{MacWilliams} the set of $Q$-tensions has Hamming weight enumerator
$${\rm hwe}(\mbox{\rm
  $Q$-tensions of $G$};t)=(t-1)^{r(G)}T\left(G;\frac{t+q-1}{t-1},t\right).$$

\section{The Tutte polynomial on $H_2$ and $H_4$}\label{H2 H4}

From equation \eqref{qflows tutte}, the Tutte polynomial specialises on $H_2$ to the Hamming weight
enumerator of the space $\mathcal{C}_2$ of $\mathbb{F}_2$-flows,  
$${\rm
  hwe}(\mathcal{C}_2;t)=(t-1)^{n(G)}T\left(G;t,\frac{t+1}{t-1}\right).$$
If $G$ is eulerian then the all $1$ vector belongs to $\mathcal{C}_2$
with the consequence that ${\rm hwe}(\mathcal{C}_2;t)=t^{|E|}{\rm
  hwe}(\mathcal{C}_2;t^{-1})$, whence if $G$ is eulerian then
\be\label{eul reciprocal}T\left(G;t,\frac{t+1}{t-1}\right)=(-1)^{n(G)}t^{r(G)}T\left(G;\frac{1}{t},\frac{1+t}{1-t}\right).\ee
Dually, if $G$ is bipartite then the all $1$ vector belongs to
$\mathcal{C}_2^\perp$ and in this case
 $$T\left(G;\frac{t+1}{t-1},t\right)=(-1)^{r(G)}t^{n(G)}T\left(G;\frac{1+t}{1-t},\frac{1}{t}\right).$$

Using the MacWilliams duality theorem~\eqref{MacWilliams},    
$$(t-1)^{n(G)}T\left(G;t,\frac{t+1}{t-1}\right)=\sum_{\mbox{\tiny
    eulerian}\,A\subseteq E}t^{|E|-|A|}=2^{-r(G)}(t-1)^{|E|}\sum_{\mbox{\tiny cutsets}\, A\subseteq E}\left(\frac{t+1}{t-1}\right)^{|E|-|A|}.$$
In particular
$$\sum_{\mbox{\tiny eulerian}\,A\subseteq E}(-1)^{|A|}=2^{|E|-|V|}P(G;2).$$

Likewise, the Tutte polynomial on $H_4$ is the Hamming weight
enumerator of the space $\mathcal{C}_4\cong\mathcal{C}_2\times\mathcal{C}_2$
 of $\mathbb{F}_4$-flows,
$${\rm  hwe}(\mathcal{C}_2\times\mathcal{C}_2;t)=(t-1)^{n(G)}T\left(G;t,\frac{t+3}{t-1}\right).$$
The MacWilliams duality theorem~\eqref{MacWilliams} here is that
$${\rm  hwe}(\mathcal{C}_2\times\mathcal{C}_2;t)=\frac{(t-1)^{|E|}}{|\mathcal{C}_2^{\perp}|^2}{\rm
  hwe}\left(\mathcal{C}_2^\perp\times\mathcal{C}_2^\perp;\frac{t+3}{t-1}\right),$$
which in terms of eulerian subgraphs of $G$ says that
$$\sum_{\mbox{\tiny
    eulerian}\,A,B\subseteq E}t^{|E|-|A\cup B|}=4^{-r(G)}(t-1)^{|E|}\sum_{\mbox{\tiny
    cutsets}\, A, B\subseteq
  E}\left(\frac{t+3}{t-1}\right)^{|E|-|A\cup B|}.$$
In particular
$$\sum_{\mbox{\tiny
    eulerian}\,A,B\subseteq E}(-3)^{|E|-|A\cup B|}=(-1)^{|E|}4^{|E|-|V|}P(G;4).$$

\subsection{Coset weight enumerators}
In a previous article~\cite{CCC07} the following specialisations of
the Tutte polynomial to the hyperbola $H_2$ and
the hyperbola $H_4$ were derived as an illustration of the techniques
afforded by elementary Fourier analysis. See also \cite{E98} (quoted in
\cite{Mi05}) for these identities in the context of coding theory. In
this section we give a combinatorial interpretation of
these identities for particular values of $t$ and derive evaluations
of the Tutte polynomial on $H_2$ and $H_4$.

When writing $\mathcal{C}_2+z\in\mathbb{F}_2^E/\mathcal{C}_2$ in the range
  of summations below, we assume that $z\in\mathbb{F}_2^E$ ranges
  over a transversal of the cosets, each coset $\mathcal{C}_2+z$
  appearing only once.

\begin{lem}\label{sum cubes} 
 Let $G=(V,E)$ be a graph and let $\mathcal{C}_2$ be the subspace of
  $\mathbb{F}_2$-flows of $G$. Then, for $t\in\mathbb{C}$,
$$\sum_{\mathcal{C}_2+z\in \mathbb{F}_2^E/\mathcal{C}_2}|{\rm
  hwe}(\mathcal{C}_2+z;t)|^2=(t+\ol{t})^{r(G)}|t-1|^{2n(G)}T\left(\!G;\frac{|t|^2+1}{t+\ol{t}},\left|\frac{t+1}{t-1}\right|^2\right),$$
$$\sum_{\mathcal{C}_2+z\in \mathbb{F}_2^E/\mathcal{C}_2}{\rm
  hwe}(\mathcal{C}_2+z;t)^2=(2t)^{r(G)}(t-1)^{2n(G)}T\left(\!G;\frac{t^2+1}{2t},\left(\frac{t+1}{t-1}\right)^2\right),$$
and 
$$\sum_{\mathcal{C}_2+z\in \mathbb{F}_2^E/\mathcal{C}_2}{\rm hwe}(\mathcal{C}_2+z;t)^3=(t+1)^{|E|} t^{r(G)}(t-1)^{2n(G)}T\left(\!G;\frac{t^2-t+1}{t},\left(\frac{t+1}{t-1}\right)^2\right).$$
\end{lem}

Note that $\left(\frac{t+1}{t-1}\right)^2$ is real if and only if
  $t\in\mathbb{R}$ or $|t|=1$, since
$\frac{t+1}{t-1}=\frac{\ol{t}-t+|t|^2-1}{|t|^2+1-t-\ol{t}}$.
By putting $t=e^{i\theta}$ in the identities of Lemma~\ref{sum cubes}, routine
calculations yield the following.\footnote{From~\cite[Corollary 7.5]{CCC07} it is readily seen
that for $r\in\mathbb{N}$ the sum of $r$th powers of the coset weight enumerators of $\mathcal{C}_2$ is a specialisation of the {\em complete weight
enumerator} (see section~\ref{poset} for a definition) of the $(r-1)$-fold direct product
$\mathcal{C}_2^\perp\times\cdots\times\mathcal{C}_2^\perp$ (isomorphic to the space of $\mathbb{F}_{2^{r-1}}$-tensions of
$G$).
Only for
$r\leq 3$ is this specialisation a Hamming weight enumerator,
and so the sum of $r$th powers of coset Hamming weight enumerators of
$\mathcal{C}_2$ is a specialisation of the Tutte polynomial only for
$r\leq 3$. For example, the sum of fourth powers of the coset weight
enumerators ${\rm hwe}(\mathcal{C}_2+z;t)$ turns out to be equal to 
$$2^{-3r(G)}(t^2-1)^{2|E|}\sum_{\mbox{\rm \tiny eulerian }\, A,B,C\subseteq E}s^{|E|-|A\cup
  B\cup C|-|A\cap B\cap C|}$$
where $s=\left(\frac{t+1}{t-1}\right)^2$.}

\begin{cor}\label{real sums} If $t=e^{i\theta}$ for some
  $\theta\in(0,2\pi)$ then
\be\label{abs sum squares theta}2^{-|E|}\sum_{\mathcal{C}_2+z\in\mathbb{F}_2^E}|{\rm hwe}(\mathcal{C}_2+z;e^{i\theta})|^2=(\cos\theta)^{r(G)}(1-\cos\theta)^{n(G)}T\left(\!G;\frac{1}{\cos\theta},\frac{1+\cos\theta}{1-\cos\theta}\right),\ee
\be\label{sum squares real}(2e^{i\theta})^{-|E|}\,\sum_{\mathcal{C}_2+z\in\mathbb{F}_2^E/\mathcal{C}_2}{\rm
  hwe}(\mathcal{C}_2+z;e^{i\theta})^2=(\cos\theta-1)^{n(G)}T\left(\!G;\cos
  \theta,\frac{\cos\theta+1}{\cos\theta-1}\right),\ee
and
$$\hspace{-5cm} (2e^{i\theta})^{-\frac{3}{2}|E|}\,\sum_{\mathcal{C}_2+z\in\mathbb{F}_2^E/\mathcal{C}_2}{\rm
  hwe}(\mathcal{C}_2+z;e^{i\theta})^3=$$
\be\label{sum cubes real}2^{-r(G)}(1+\cos\theta)^{\frac{1}{2}|E|}(\cos\theta-1)^{n(G)}T\left(\!G;2\cos \theta-1,\frac{\cos\theta+1}{\cos\theta-1}\right).\ee
\end{cor}


If $\theta\in\mathbb{Q}$ then $2 \cos \theta =e^{i\theta}+e^{-i\theta}$ is an algebraic
integer, and hence an ordinary integer if $\cos\theta\in\mathbb{Q}$. 
The only rational values of $\cos \theta$ for $\theta\in\mathbb{Q}$
are thus
$0,\pm\frac{1}{2},\pm 1$, corresponding to $\theta=\pm\frac{\pi}{2},
\pm\frac{\pi}{3}, \pm\frac{2\pi}{3},\pi,0$.
Thus when $e^{i\theta}$ is a $q$th root of unity for $q\in\{2,3,4,6\}$ the evaluations
of the Tutte polynomial in Corollary~\ref{real sums} are at rational
points (see Table~\ref{table} below). 

\begin{table}[ht]
\begin{center}
\caption{{\bf \small Evaluations of the Tutte polynomial in the
  identities of Corollary~\ref{real sums}}}\label{table}
$$\begin{array}{c|c|c|c|c}
\theta & \cos \theta & \mbox{\rm eq. \eqref{abs sum squares theta}, point on $H_2$} & \mbox{\rm eq. \eqref{sum
  squares real}, point on $H_2$} &  \mbox{\rm eq. \eqref{sum cubes real}, point on $H_4$}\\
\hline
\pi & -1 & (-1,0) & (-1,0) & (-3,0)\\
2\pi/3 & -\frac{1}{2} & (-2,\frac{1}{3}) & (-\frac{1}{2},-\frac{1}{3}) &
(-2,-\frac{1}{3})\\
\pi/2 & 0 &  * & (0,-1) & (-1,-1)\\
\pi/3 & \frac{1}{2} & (2,3) & (\frac{1}{2},-3) & (0,-3)\end{array}$$
{\small * For $\theta=\pi/2$, the right-hand side of equation (\ref{abs sum
  squares theta}) is equal to $1$ independent of $G$.

For $\theta=\pi$ the factor
$(1+\cos\theta)^{\frac{1}{2}|E|}$ on the right-hand side of (\ref{sum cubes real}) is equal
to zero.}

\end{center}
\end{table}

In order to give combinatorial interpretations of identity~\eqref{abs sum
  squares theta} (and identity~\eqref{sum squares real}) we shall be interested in the correlation between two types of
event when choosing $A,B\subseteq E$ uniformly at random. First, the event that $A\bigtriangleup B$ is eulerian. Second,
for various integers $q$, 
the event that $|A|-|B|$ (respectively $|A|+|B|$) belongs to a certain subset of congruence
classes modulo $q$. 

Similarly, in order to interpret identity~\eqref{sum cubes real} we 
choose $A,B,C\subseteq E$ uniformly at random and look at the correlation between
the event that $A\bigtriangleup B, B\bigtriangleup C,
C\bigtriangleup A$ are each
eulerian and the event that $|A|+|B|+|C|$ belongs to a certain
subset of congruence classes modulo $q$.

\subsection{Bias}

If $q$ is not a power of two, for fixed $k\in\{0,1,\ldots, q-1\}$
none of the events $|A|-|B|\equiv k (\bmod\, q)$, $|A|+|B|\equiv k (\bmod\, q)$ or $|A|+|B|+|C|\equiv
  k (\bmod\, q)$ can have
  probability $\frac{1}{q}$ when $A,B,C\subseteq E$ are taken
  uniformly at random. As observed in Remark~\ref{equidistributed}
  below, amongst powers of two only for $q=2$ is it true
  that the values of $|A\pm|B|(+|C|)$ are equidistributed modulo $q$ (although it remains possible that for some values
 of $k$ the event $|A\pm|B|(+|C|)\equiv k (\bmod\, q)$ has probability $\frac1q$).

Let $\Sigma$ be an event in the uniform probability space on pairs
$A,B\subseteq E$ or triples $A,B,C\subseteq E$, and $\ol{\Sigma}$ is its
complement. In the sequel the event $\Sigma$ takes the form $|A|\pm
|B|\in S (\bmod\, q)$ or
$|A|\!+\!|B|\!+\!|C|\in S (\bmod\, q)$ for a subset $S$ of the integers $\{0,1\ldots, q-1\}$ modulo
$q$. 

Define the {\em bias} towards $\Sigma$ by
$${\rm Bias}(\Sigma)=\mathbb{P}(\Sigma)-\mathbb{P}(\ol{\Sigma}).$$
Note that ${\rm Bias}(\ol{\Sigma})=-{\rm Bias}(\Sigma)$ and the event $\Sigma$ has probability
$\frac{1}{2}[1+{\rm Bias}(\Sigma)]$.

For example, if $q=2$ and $S=\{0\}$ then 
$\Sigma=\{A,B\subseteq E:|A|+|B|\in S (\bmod\, q)\}$ is the event that $|A|+|B|$
is even and this has the same probability as the event
$\ol{\Sigma}=\{A,B\subseteq E:|A|+|B|\not\in S (\bmod\, q)\}$ that
$|A|+|B|$ is odd. For $q>2$, the bias of an event of the form $|A|\pm
|B|\in S (\bmod\, q)$ or
$|A|+|B|+|C|\in S (\bmod\, q)$ is usually not equal to zero. This is excepting the case
when $q$ is even and $S=\{0,2,\ldots, q-2\}$ or $S=\{1,3,\ldots, q-1\}$ for
 which the event $\Sigma$ is about the parity of $|A|\pm|B|(+|C|)$
 again. However, in Theorem~\ref{A + B + C with q=4} it is found that 
${\rm Bias}(|A|\!+\!|B|\!+\!|C|\equiv 0,1 (\bmod\, 4))=0$ when
 $|E|\equiv 1 (\bmod\, 4)$.

\begin{lem}\label{bias S}
Suppose $A,B,C\subseteq E$ are chosen uniformly at
  random and that $S\subseteq\{0,1\ldots, q-1\}$. Then 
$${\rm Bias}(\Sigma)=2q^{-1}\langle\wh{g},\wh{1_S}\rangle-1,$$
where $1_S$ is the indicator function of $S$ and
$$\wh{g}(k)=\begin{cases} 2^{-|E|}(1+\cos\frac{2\pi k}{q})^{|E|}\\
 2^{-|E|}e^{-2\pi ik|E|/q}(1+\cos\frac{2\pi k}{q})^{|E|} \\
 2^{-\frac{3}{2}|E|}e^{-2\pi ik\frac{3}{2}|E|/q}(1+\cos\frac{2\pi k}{q})^{\frac{3}{2}|E|}\end{cases}$$
according as
$$\Sigma=\begin{cases}\{A,B\subseteq E:|A|-|B|\in S(\bmod\, q)\}\\
\{A,B\subseteq E:|A|+|B|\in S(\bmod\, q)\}\\
\{A,B,C\subseteq E:|A|+|B|+|C|\in S(\bmod\,q)\}\end{cases}.$$
\end{lem}

\begin{pf} We prove the lemma for the event
  $\Sigma=\{A,B\subseteq E:\,|A|-|B|\in S(\bmod\, q)\}$. The
  other cases are similar. 

Define
$$g(\ell)=\mathbb{P}\big(\,|A|-|B|\equiv\ell(\bmod\, q)\,\big).$$
Then 
$$\mathbb{P}(\Sigma)=\sum_{\ell\in
  S}g(\ell)=\langle g,1_S\rangle=q^{-1}\langle\wh{g},\wh{1_S}\rangle,$$
using Plancherel's formula at the end.
By definition, ${\rm
  Bias}(\Sigma)=2\mathbb{P}(\Sigma)-1$, and the first part of
the lemma is proved. It remains to calculate $\wh{g}(k)$ for $k\in\{0,1,\ldots, q-1\}$:
\begin{align*}
\wh{g}(k)=\sum_{0\leq \ell\leq q-1}e^{-2\pi i k\ell/q}g(\ell) & =2^{-2|E|}\sum_{A, B\subseteq E}e^{-2\pi ik(|A|-|B|)/q}\\
 & =2^{-2|E|}\sum_{A\subseteq
  E}e^{-2\pi ik|A|/q}\sum_{B\subseteq
  E}e^{2\pi ik|B|/q}\\
 & = 2^{-2|E|}(1+e^{-2\pi ik/q})^{|E|}(1+e^{2\pi ik/q})^{|E|}\\
 & =2^{-2|E|}(e^{2\pi ik/q}+2+e^{-2\pi ik/q})^{|E|}\\
 & =2^{-|E|}(1+\cos\frac{2\pi k}{q})^{|E|}.\end{align*}
This completes the proof.
\end{pf}

\begin{rem}\label{equidistributed}
 In the notation of the proof of Lemma~\ref{bias S}, the probabilities
$g(\ell)=\mathbb{P}(|A|-|B|=\ell(\bmod\, q))$ are equal for all
 $\ell\in\{0,1,\ldots, q-1\}$ if and only if 
$$\sum_{0\leq\ell\leq q-1}g(\ell)e^{-2\pi i\ell/q}=0.$$
But the left-hand sum is equal to
$\wh{g}(1)=2^{-|E|}(1+\cos\frac{2\pi}{q})^{|E|}$, and this is equal to
zero if and only if $q=2$. Hence the values $|A|\!-\!|B|$ for $A,B\subseteq E$ are
equidistributed modulo $q$ if and only if $q=2$. Similarly, the values
$|A|\!+\!|B|$ for $A,B\subseteq E$ and the values of
$|A|\!+\!|B|\!+\!|C|$ for $A,B,C\subseteq E$ are only equidistributed
modulo $q$ when $q=2$.
\end{rem}

Let $\Delta$ be an event in the uniform probability space on pairs
$A,B\subseteq E$ or triples $A,B,C\subseteq E$.
We define the
conditional bias of $\Sigma$ given $\Delta$ by 
$${\rm Bias}(\Sigma\,|\,\Delta)=\mathbb{P}(\Sigma\,|\,\Delta)-\mathbb{P}(\ol{\Sigma}\,|\,\Delta).$$
In what follows, $\Delta$ is either the event that $A\bigtriangleup B$ is
eulerian (where $\Sigma$ is one of the events $|A|\pm|B|\in
S\,(\mbox{\rm mod}\, q)$ for
some $S\subseteq\{0,1\ldots, q-1\}$) or the event
that $A\bigtriangleup B$ and $B\bigtriangleup C$ are both eulerian
(where $\Sigma$ is the event that $|A|+|B|+|C|\in S\,(\mbox{\rm mod}\,
q)$ for some $S\subseteq\{0,1\ldots, q-1\}$).

The covariance of (the indicator functions of) the events $\Sigma$ and
$\Delta$ is defined by the difference
$\mathbb{P}(\Sigma\cap\Delta)-\mathbb{P}(\Sigma)\mathbb{P}(\Delta)$.
We define the {\em correlation}\footnote{The {\em correlation
    coefficient} of the indicator functions of the events  $\Sigma$
  and $\Delta$ is another normalisation of the
  covariance, namely $$\frac{\mathbb{P}(\Sigma\cap
    \Delta)-\mathbb{P}(\Sigma)\mathbb{P}(\Delta)}{\sqrt{\mathbb{P}(\Sigma)\mathbb{P}(\Delta)(1-\mathbb{P}(\Sigma))(1-\mathbb{P}(\Delta))}}.$$}
between the events $\Sigma$
and $\Delta$ by dividing the covariance through by
$\mathbb{P}(\Delta)$,
$${\rm Correlation}(\Sigma\,|\,\Delta)=\mathbb{P}(\Sigma\,|\,\Delta)-\mathbb{P}(\Sigma).$$

Correlation is related to bias via the relation
$${\rm Bias}(\Sigma\,|\,\Delta)-{\rm Bias}(\Sigma)=2\,{\rm
  Correlation}(\Sigma\,|\,\Delta).$$
When ${\rm Bias}(\Sigma)\neq 0$ we shall have occasion to also measure correlation via the ratio
\be\label{ratio corr}\frac{{\rm Bias}(\Sigma\,|\,\Delta)}{{\rm Bias}(\Sigma)}=\frac{2{\rm
    Correlation}(\Sigma\,|\,\Delta)}{{\rm Bias}(\Sigma)}+1.\ee
This can be viewed as the scale factor from the existing
bias towards $\Sigma$ to the bias towards $\Sigma$ given the event $\Delta$.
When the ratio~\eqref{ratio corr} is greater than $1$ the existing bias towards $\Sigma$ is magnified,
when the ratio is less than $1$ the existing bias is diminished.
The ratio~\eqref{ratio corr} is equal to
$1$ if and only if there is no correlation between the events
$\Sigma$ and $\Delta$. 

In the next section we derive general expressions for ${\rm
  Bias}(\Sigma)$ and ${\rm
  Bias}(\Sigma\,|\,\Delta)$, for any choice of $q$ and $S$ in the definition of $\Sigma$,
the latter expressed in terms of the Tutte
polynomial evaluations of Corollary~\ref{real sums} for
$\theta\in\{2\pi k/q:k=1,\ldots, q-1\}$.
By taking $q\in\{2,3,4,6\}$ we obtain interpretations for various
evaluations of the Tutte polynomial at the points listed in Table~\ref{table}. The reason why
we limit ourselves to $q\in\{2,3,4,6\}$ is not specifically on account of these
corresponding to evaluations at rational points, but rather that only for
these values of $q$ do we obtain evaluations of the Tutte
polynomial at a single point rather than a sum of evaluations at two
or more distinct points.

An evaluation of the Tutte
polynomial at a single point is of interest as such an evaluation is a {\em
  Tutte-Grothendieck invariant} on graphs, i.e., satisfies a ``linear''
deletion-contraction recurrence relation. 
Specifically, if $\setminus$ denotes deletion and $/$ contraction, a function $f$ on graphs satisfying
$$f(G)=\begin{cases}
 af(G\setminus e)+bf(G/e) &  \mbox{ if $e$ is neither an isthmus nor a loop,}\\
 xf(G\setminus e) & \mbox{ if $e$ is an isthmus,}\\
 yf(G\setminus e) &    \mbox{ if $e$ is a loop,}\end{cases}$$
and with value $c^{|V|}$ on the edgeless graph $(V,\emptyset)$ is
given by the evaluation\footnote{An {\em isthmus} or bridge is an edge
  forming a cutset of size $1$ and a {\em loop} is an edge forming a
   cycle of size $1$. If $b=0$ then 
$$f(G)=c^{|V|}a^{n(G)-\ell(G)}x^{r(G)}y^{\ell(G)},$$
and if $a=0$ then $$f(G)=c^{k(G)+i(G)}b^{r(G)-i(G)}x^{i(G)}y^{n(G)},$$
where $i(G)$ is the number of isthmuses and $\ell(G)$ the number of
loops in $G$.}  \cite{B98}, \cite{DW93},
$$f(G)=c^{k(G)}a^{n(G)}b^{r(G)}T(G;cx/b, y/a).$$
\label{fn}
Many of the evaluations of the Tutte polynomial in sections~\ref{square root of unity} to \ref{sixth root of unity} have
been highlighted as theorems either because they have other well-known combinatorial meanings (such
as the number of nowhere-zero $4$-flows in Theorem~\ref{sum cubes with
  q=6}) and the opaqueness of the connection to these other
interpretations is intriguing,
or because they on the contrary do not have such other well-known interpretations
(such as $T(G;-2,-\frac{1}{3})$ in Theorem~\ref{cube root A
  B C}).
Direct proofs of the corresponding
deletion-contraction recurrences do not seem straightforward in many cases.

\subsection{Evaluations of coset weight enumerators at $q$th roots of
  unity}\label{evaluation at qth}

In this section the identities in
Corollary~\ref{real sums} are given interpretations in terms of ${\rm
  Bias}(\Sigma|\Delta)$ where $\Sigma$ takes the form
$|A|\pm|B|(+|C|)\in S\,(\mbox{\rm mod}\,q)$ and $\Delta$ is the event
that $A\bigtriangleup B$ (and $B\bigtriangleup C$) is eulerian. These interpretations in their
general form make for rather dull reading, but their particular cases for $q\in\{2,3,4,6\}$ are the more
interesting theorems that follow as corollaries. The latter are
presented in sections~\ref{square root of unity} to
\ref{sixth root of unity} to which the reader might wish to turn
before referring back to this section.

\begin{lem}\label{interpret h2|} Suppose $A,B\subseteq E$ are chosen uniformly at
  random and $\Delta$ is the event that $A\bigtriangleup B$ is
  eulerian. Then, for any $k\in\{0,1,\ldots, q-1\}$,
$$2^{-|E|}\sum_{\mathcal{C}_2+z\in\mathbb{F}_2^E/\mathcal{C}_2}|{\rm
  hwe}(\mathcal{C}_2+z;e^{2\pi ik/q})|^2=2^{n(G)}\sum_{0\leq\ell\leq
  q-1}e^{2\pi i k\ell/q}\,\mathbb{P}(|B|-|A|\equiv \ell(\bmod{q})\,\mid\,\Delta).$$
\end{lem}
\begin{pf}
Given a coset $\mathcal{C}_2+z$ and $\ell\in\{0,\ldots,q-1\}$, define
 $$p_\ell=p_\ell(\mathcal{C}_2+z)=\#\{x\in\mathcal{C}_2+z:|E|-|x|\equiv
\ell(\bmod q)\},$$
so that
$$\sum_{0\leq
  \ell\leq q-1}e^{2\pi ik\ell/q}p_\ell={\rm hwe}(\mathcal{C}_2+z;e^{2\pi ik/q}).$$
Going on to define
$$P_\ell=P_\ell(\mathcal{C}_2+z)=\sum_{j-k\equiv \ell(\bmod
  q)}p_jp_k,$$
we have
\be\label{|hwe|^2 to Pl}|{\rm hwe}(\mathcal{C}_2+z;e^{2\pi ik/q})|^2 = \big|\sum_{0\leq
  \ell\leq q-1}e^{2\pi ik\ell/q}p_\ell\big|^2= \sum_{0\leq\ell\leq q-1} e^{2\pi ik\ell/q}P_\ell.\ee

Let $C\subseteq E$ have indicator vector $z\in\mathbb{F}_2^E$ and suppose $A\subseteq
E$ is chosen uniformly at random with indicator vector $x\in\mathbb{F}_2^E$.
Then $|E|-|x|=|E\setminus A|$ and $x\in\mathcal{C}_2+z$ if and only if
$x+z\in\mathcal{C}_2$, i.e., $A\bigtriangleup C$ is eulerian. Hence
$$2^{-|E|}p_\ell(\mathcal{C}_2+z)=\mathbb{P}(|E\setminus
A|\equiv\ell(\bmod\, q)\cap A\bigtriangleup C\,\mbox{\rm eulerian}).$$
Similarly, if $A,B\subseteq E$ are chosen uniformly at random then
$$2^{-2|E|}P_\ell(\mathcal{C}_2+z) =$$
$$\sum_{j-k\equiv \ell(\bmod q)}\mathbb{P}(|E\setminus A|\equiv j(\!\bmod q)\,\cap\,A\bigtriangleup C\,\mbox{\rm eulerian})\,\mathbb{P}(|E\setminus B|\equiv k(\!\bmod
q)\,\cap B\bigtriangleup C\,\mbox{\rm eulerian})$$
$$ = \mathbb{P}(|E\setminus A|-|E\setminus B|\equiv\ell(\bmod q)\,\cap
A\bigtriangleup C\,\mbox{\rm eulerian}\cap B\bigtriangleup
C\,\mbox{\rm eulerian}).$$

Given $C,C'\subseteq E$, the events $\{A\bigtriangleup C\,\mbox{\rm
  eulerian}\cap B\bigtriangleup C\,\mbox{\rm eulerian}\}$ and $\{A\bigtriangleup C'\,\mbox{\rm
  eulerian}\cap B\bigtriangleup C'\,\mbox{\rm eulerian}\}$ are
either equal (when $C\bigtriangleup C'$ is eulerian) or
disjoint. For suppose that $C$ has indicator vector
$z$ and $C'$ indicator vector $z'$. If
$x+z\in\mathcal{C}_2$ and $x+z'\in\mathcal{C}_2$ then
$z+z'\in\mathcal{C}_2$, i.e., $C\bigtriangleup C'$ is eulerian. Conversely, if $z+z'\in\mathcal{C}_2$ and $x+z\in\mathcal{C}_2$ then $x+z+(z+z')=x+z'\in\mathcal{C}_2$.

Letting $C\subseteq E$ range over a collection of subsets no
two of which have eulerian symmetric difference, the union of events $\{A\bigtriangleup C\,\mbox{\rm
  eulerian}\cap B\bigtriangleup C\,\mbox{\rm eulerian}\}$ is thus a
disjoint union and  equal
to the event $\{A\bigtriangleup B\,\mbox{\rm
  eulerian}\}=\Delta$. Hence, letting $z$ range over a transversal of cosets of $\mathcal{C}_2$,
$$2^{-2|E|}\sum_{\mathcal{C}_2+z\in\mathbb{F}_2^E/\mathcal{C}_2}P_\ell(\mathcal{C}_2+z)=\mathbb{P}(|E\setminus
A|-|E\setminus B|\equiv \ell(\bmod\, q)\cap\Delta).$$
From equation~\eqref{|hwe|^2 to Pl} it follows that
$$2^{-2|E|}\sum_{\mathcal{C}_2+z\in\mathbb{F}_2^E/\mathcal{C}_2}|{\rm hwe}(\mathcal{C}_2+z;e^{2\pi ik/q})|^2=\sum_{0\leq\ell\leq
  q-1}e^{2\pi i k\ell/q}\,\mathbb{P}(|B|-|A|\equiv\ell\,(\bmod\,q)\,\cap\,\Delta).$$
Dividing through by $\mathbb{P}(\Delta)=2^{-r(G)}$ gives the result.
\end{pf}

\begin{lem}\label{interpret h2} Suppose that $A,B\subseteq E$ are chosen uniformly at
  random and $\Delta$ is the event that $A\bigtriangleup B$ is eulerian. Then, for any $k\in\{0,1,\ldots, q-1\}$,
$$2^{-|E|}\sum_{\mathcal{C}_2+z\in\mathbb{F}_2^E/\mathcal{C}_2}{\rm
  hwe}(\mathcal{C}_2+z;e^{2\pi ik/q})^2=2^{n(G)}\sum_{0\leq\ell\leq
  q-1}e^{2\pi i k\ell/q}\,\mathbb{P}(|A|+|B|\equiv \ell\,(\mbox{\rm
  mod}\, q)\,\mid\,\Delta).$$
\end{lem}
\begin{pf}
The same mutatis mutandis as the proof of Lemma~\ref{interpret
  h2|} with 
$$P_\ell=P_\ell(\mathcal{C}_2+z)=\sum_{j+k\equiv\ell(\bmod q)}p_jp_k$$
replacing the definition of $P_\ell$ given in the proof of that lemma.
 Note that since $A,B\subseteq E$ are chosen uniformly at
random and $(E\setminus A)\bigtriangleup(E\setminus B)=A\bigtriangleup
  B$, by symmetry we have $\mathbb{P}(|E\setminus
A|+|E\setminus B|\equiv\ell\,(\mbox{\rm mod}\,
  q)\,\cap\,\Delta)=\mathbb{P}(|A|+|B|\equiv\ell\,(\mbox{\rm mod}\,
  q)\,\cap\,\Delta)$. 
\end{pf}

\begin{lem}\label{interpret h3} Suppose that $A,B,C\subseteq E$ are chosen uniformly at
  random and $\Delta$ is the event that $A\bigtriangleup B,
  B\bigtriangleup C$ are both eulerian. Then, for any $k\in\{0,1,\ldots, q-1\}$,
$$2^{-|E|}\sum_{\mathcal{C}_2+z\in\mathbb{F}_2^E/\mathcal{C}_2}{\rm
  hwe}(\mathcal{C}_2+z;e^{2\pi ik/q})^3=4^{n(G)}\sum_{0\leq\ell\leq
  q-1}e^{2\pi i k\ell/q}\,\mathbb{P}(|A|+|B|+|C|\equiv \ell\,(\mbox{\rm
  mod}\, q)\,\mid\,\Delta).$$
\end{lem}
\begin{pf} The same as the proof of Lemma~\ref{interpret h2|} with
  the following being the main alterations.
Set $$P_\ell=P_\ell(\mathcal{C}_2+z)=\sum_{i+j+k\equiv \ell\,(\mbox{\tiny mod}\,
  q)}p_ip_jp_k.$$
Then
$${\rm hwe}(\mathcal{C}_2+z;e^{2\pi i k/q})^3=\sum_{0\leq\ell\leq
  q-1} e^{2\pi ik\ell/q}P_\ell$$
and
$$2^{-3|E|}\sum_{\mathcal{C}_2+z\in\mathbb{F}_2^E/\mathcal{C}_2}P_{\ell}(\mathcal{C}_2+z)=\mathbb{P}(|E\setminus
A|+|E\setminus B|+|E\setminus C|\equiv\ell\,(\bmod\,q)\,\cap\,\Delta).$$
As in the proof of Lemma~\ref{interpret h2}, we can by symmetry replace $E\setminus A, E\setminus B, E\setminus C$
by $A,B,C$. Hence
$$2^{-3|E|}\sum_{\mathcal{C}_2+z\in\mathbb{F}_2^E/\mathcal{C}_2}{\rm
  hwe}(\mathcal{C}_2+z;e^{2\pi ik/q})^3=\sum_{0\leq\ell\leq
  q-1}e^{2\pi i k\ell/q}\,\mathbb{P}(|A|+|B|+|C|\equiv\ell\,(\bmod\,q)\,\cap\,\Delta).$$
Dividing through by $\mathbb{P}(\Delta)=2^{-2r(G)}$ gives the result.
\end{pf}

\begin{lem}\label{Bias Sigma Delta} Suppose that $A,B\,(,C)\subseteq E$ are chosen uniformly at
  random and $\Delta$ is the event that $A\bigtriangleup B$ (and
  $B\bigtriangleup C$) is eulerian. Then
$${\rm Bias}(\Sigma\,\mid\,\Delta)=2q^{-1}\langle\wh{f},\wh{1_S}\rangle-1,$$
where $1_S$ is the indicator function of $S$ and, for
$k\in\{0,1,\ldots, q-1\}$,
$$\wh{f}(k)=\begin{cases} 2^{-n(G)}\left(\cos\frac{2\pi
      k}{q}\right)^{r(G)}\left(1\!-\!\cos\frac{2\pi
      k}{q}\right)^{n(G)}T\left(\!G;\frac{1}{\cos\frac{2\pi
        k}{q}},\frac{1+\cos\frac{2\pi k}{q}}{1-\cos\frac{2\pi k}{q}}\right)\\
 2^{-n(G)}e^{-2\pi i k|E|/q}\left(\cos\frac{2\pi
     k}{q}\!-\!1\right)^{n(G)}T\left(\!G;\cos\frac{2\pi
     k}{q},\frac{\cos\frac{2\pi k}{q}+1}{\cos\frac{2\pi
       k}{q}-1}\right)\\
2^{-n(G)-\frac{1}{2}|E|}e^{-2\pi i k\frac{3}{2}|E|/q}\left(1\!+\!\cos\frac{2\pi
     k}{q}\right)^{\frac{1}{2}|E|}\left(\cos\frac{2\pi
     k}{q}\!-\!1\right)^{n(G)}T\left(\!G;2\cos\frac{2\pi
     k}{q}\!-\!1,\frac{\cos\frac{2\pi k}{q}+1}{\cos\frac{2\pi
       k}{q}-1}\right)\end{cases}$$
according as
 $$\Sigma=\begin{cases}\{A,B\subseteq E:|A|-|B|\in S\,(\bmod\,q)\}\\
\{A,B\subseteq E:|A|+|B|\in S\,(\bmod\,q)\}\\
\{A,B,C\subseteq E: |A|\!+\!|B|\!+\!|C|\in S\,(\bmod\,q)\}\end{cases}.$$
\end{lem}

\begin{pf}  Again we prove the result for $\Sigma=\{A,B\subseteq
  E:|A|-|B|\in S\,(\bmod\,q)\}$ and $\Delta$ the event that
  $A\bigtriangleup B$ is eulerian, the other cases being entirely
  similar.

Define $f(\ell)=\mathbb{P}(|A|-|B|\equiv \ell\,\mid\,\Delta)$. Then
$\mathbb{P}(\Sigma)=\langle
f,1_S\rangle=q^{-1}\langle\wh{f},\wh{1_S}\rangle$, by Parseval's
formula. By Lemma~\ref{interpret h2|}, 
$$\wh{f}(k)=2^{-|E|-n(G)}\sum_{\mathcal{C}_2+z\in\mathbb{F}_2^E/\mathcal{C}_2}|{\rm
  hwe}(\mathcal{C}_2+z;e^{-2\pi i k/q})|^2$$
and identity~\eqref{abs sum squares theta} gives the result. 
\end{pf}

Lemma~\ref{Bias Sigma Delta} shows that ${\rm
  Bias}(\Sigma\,\mid\,\Delta)$ is given in terms of evaluations of the Tutte
polynomial at one or more points.  The remainder of this section is
  spent establishing when an evaluation at just one point is involved
  and a Tutte-Grothendieck invariant results. 

 Note that $\wh{f}(0)=1$ for each $\wh{f}$ defined in Lemma~\ref{Bias
  Sigma Delta} ($f$ defines a probability distribution on
  $\mathbb{Z}_q$ and $\wh{f}(0)=\sum f(k)=1$). Recall the definition
  of $g$ from
  Lemma~\ref{bias S}, where Bias$(\Sigma)$ is expressed in terms of the
  inner product $\langle\wh{g},\wh{1_S}\rangle$. 

 Since $f$ and $g$ are
  real-valued, $\wh{g}(-k)=\ol{\wh{g}(k)}$ and $\wh{f}(-k)=\ol{\wh{f}(k)}$. Remark
  also that if $q$ is even and $\Sigma=\{A,B,C\subseteq
  E:|A|\!+\!|B|\!+\!|C|\in S\,(\bmod\,q)\}$ then
  $\wh{f}(q/2)=0$ (for all graphs $G$), but that $\wh{f}(k)\neq 0$ for
  some graph $G$ when $k\neq q/2$.
 
The {\em support} of a function $h:Q\rightarrow\mathbb{C}$ is defined
by ${\rm supp}(h)=\#\{k\in Q:h(k)\neq 0\}$. Thus ${\rm
  supp}(\wh{g})=\mathbb{Z}_q\setminus\{q/2\}$, and ${\rm
  supp}(\wh{f})=\mathbb{Z}_q$ (or possibly
$\mathbb{Z}_q\setminus\{q/2\}$ as we have just seen).
From Lemma~\ref{Bias Sigma Delta}, ${\rm
  Bias}(\Sigma\,|\,\Delta)=2q^{-1}\langle\wh{f},\wh{1_S}\rangle-1$ will involve an evaluation of the Tutte
  polynomial at a single point (valid for all graphs $G$) only if
  ${\rm supp}(\wh{f}\cdot \wh{1_S})\subseteq\{0,\ell,-\ell\}$ for some
  $\ell$, or ${\rm supp}(\wh{f}\cdot \wh{1_S})\subseteq\{0,\ell,-\ell,q/2\}$ if $q$ is even and $\Sigma=\{A,B,C\subseteq
  E:|A|\!+\!|B|\!+\!|C|\in S\,(\bmod\,q)\}$. This is so that the
  only non-zero terms contributing to the expression
  $2q^{-1}\langle\wh{f},\wh{1_S}\rangle-1$ are $\wh{f}(\ell)\overline{\wh{1_S}(\ell)}$ and
  its complex conjugate $\wh{f}(-\ell)\overline{\wh{1_S}(-\ell)}$.

Suppose $\Sigma'$ is an event defined just as $\Sigma$ except with
  $S'\subseteq\mathbb{Z}_q\setminus S$ replacing $S$.  Then the bias towards $S$ can be compared to
  the bias toward $S'$ by considering the difference ${\rm
  Bias}(\Sigma\,|\,\Delta)-{\rm
  Bias}(\Sigma'\,|\,\Delta)$. When $S'=\mathbb{Z}_q\setminus S$
  this difference is simply $2{\rm Bias}(\Sigma\,|\,\Delta)$. The
  criterion for ${\rm Bias}(\Sigma\,|\,\Delta)-{\rm Bias}(\Sigma'\,|\,\Delta)$ to be an evaluation of the Tutte polynomial at a
  single point valid for all graphs is that ${\rm supp}(\wh{f}\cdot
  \wh{1_S-1_{S'}})\subseteq\{0,\ell,-\ell\, (,q/2)\}$ (with $q/2$
  included under the same conditions as before). 

 The multiplicative group
of units of $\mathbb{Z}_q$ is denoted by $\mathbb{Z}_q^\times$ and has
order $\phi(q)$, where $\phi(q)=\#\{1\leq k\leq q:(k,q)=1\}$ is
Euler's totient function. For $S\subseteq\mathbb{Z}_q$ and
$\ell\in\mathbb{Z}_q$ we write $\ell S=\{\ell s:s\in S\}$.
 
\begin{lem}\label{support} If $h:\mathbb{Z}_q\rightarrow\mathbb{Q}$ and $\wh{h}(k)\neq
  0$ then  ${\rm supp}(\wh{h})\supseteq k\mathbb{Z}_q^\times$. 

If ${\rm supp}(\wh{h})\subseteq d\mathbb{Z}_q$ for some divisor $d$ of $q$
  then $h$ is constant on cosets of $(q/d)\mathbb{Z}_q$.

Thus if $\wh{h}(k)\neq 0$ for a unit $k$ of $\mathbb{Z}_q$ then
$\wh{h}(\ell)\neq 0$ for all $\ell\in\mathbb{Z}_q^\times$, while if
there is no unit $k$ in the support of $\wh{h}$ then $h$ is constant on cosets of a proper subgroup of $\mathbb{Z}_q$.  
\end{lem}
\begin{pf}
The first statement depends on the fact that $h$ takes rational
values.  Suppose
$j\mapsto \sigma_j$ is the natural isomorphism of the group of Galois automorphisms of $\mathbb{Q}(e^{2\pi
  i/q})$ with $\mathbb{Z}_q^\times$, i.e.,\ $\sigma_j:e^{2\pi
  i/q}\mapsto e^{2\pi
  i j/q}$. Then
\begin{align*}\sigma_j(\wh{h}(k)) &
  =\sigma_j\big(\sum_{0\leq\ell\leq q-1}h(\ell)e^{2\pi i k\ell/q}\big)\\
 & =\sum_{0\leq \ell\leq q-1}h(\ell)\sigma_j(e^{2\pi i k\ell/q})\\
 & =\wh{h}(jk).\end{align*}
 Hence $\wh{h}(k)\neq 0$ implies $\wh{h}(jk)\neq 0$ for all $j\in\mathbb{Z}_q^\times$.

For the second statement, the assumption is that
$\wh{h}=\sum_{0\leq \ell\leq q/d-1}\wh{h}(d\ell)1_{d\ell}$. Taking Fourier transforms, for each
$k\in\mathbb{Z}_q$ we have 
$$qh(k)=\sum_{0\leq \ell\leq q/d-1}\wh{h}(d\ell)e^{2\pi i k\ell/(q/d)}.$$
Then
$$h(k+q/d)=q^{-1}\sum_{0\leq\ell\leq q/d-1}\wh{h}(d\ell)e^{2\pi i
    k\ell/(q/d)+2\pi i\ell}=h(k),$$
  so that $h$ has period $q/d$ and is constant on additive cosets of $(q/d)\mathbb{Z}_q$.
\end{pf}

Whether ${\rm  Bias}(\Sigma\,|\, \Delta)-{\rm
    Bias}(\Sigma'\,|\, \Delta)$ involves an evaluation of the Tutte
    polynomial at a single point depends on how many how zero terms
    there are in its expression as $2q^{-1}\langle\wh{f},\wh{1_S}-\wh{1_{S'}}\rangle$ obtained from
    Lemma~\ref{Bias Sigma Delta} applied to $\Sigma$ and $\Sigma'$. We require $|S|=|S'|$ for
    the sets $S, S'$ defining the events $\Sigma, \Sigma'$ since
    $\wh{1_S}(0)-\wh{1_{S'}}(0)=|S|-|S'|$ and $\wh{f}(0)=1$. 

\begin{cor}\label{when} Suppose that $\Sigma$ is one of the events $\{A,B\subseteq
    E:|A|\pm|B|\in S(\bmod q)\}$, $\Sigma'$ is similarly
    defined with $S'\subseteq\mathbb{Z}_q\setminus S$ in place of
    $S$, and $\Delta$ is the event
    that $A\bigtriangleup B$ is eulerian. Then ${\rm
    Bias}(\Sigma\,|\, \Delta)-{\rm
    Bias}(\Sigma'\,|\, \Delta)$ is up to a factor depending only on
    $|E|$ and $r(G)$ an evaluation of the Tutte
    polynomial of $G$ at a single point only if $|S|=|S'|$ and
    $q\in\{2,3,4\}$ or $S, S'$ are each unions of additive cosets of $d\mathbb{Z}_q$ for $d\in\{2,3,4\}$ a divisor
    of $q$.

If $\Sigma$ is the event $\{A,B,C\subseteq E:|A|+|B|+|C|\in S\,(\mbox{\rm mod}\,
    q)\}$, $\Sigma'$ the same event with
    $S'\subseteq\mathbb{Z}_q\setminus S$ in place of $S$,  and $\Delta$ the event
    that $A\bigtriangleup B, B\bigtriangleup C$ are both eulerian, then ${\rm
    Bias}(\Sigma\,|\, \Delta)-{\rm
    Bias}(\Sigma'\,|\, \Delta)$ involves an evaluation of the Tutte
    polynomial at a single point only if $|S|=|S'|$ and $q\in\{2,3,4,6\}$ or $S$
    and $S'$ are each unions of additive cosets of $d\mathbb{Z}_q$ for $d\in\{2,3,4,6\}$ a divisor
    of $q$.
\end{cor}
Note that if $S$ is a union of additive cosets of $d\mathbb{Z}_q$ then the
event $\Sigma$ is a congruence condition modulo $d$ so these choices for
$S$ are herewith ignored.
\begin{pf}
  Let $h=1_S-1_{S'}$. The only integers $q\geq 2$ for which $\phi(q)\leq 2$ are
  $2,3,4,6$. By Lemma~\ref{support} either we are in the case where
  $\wh{h}$ is supported on an additive subgroup $d\mathbb{Z}_q$ or
  $\wh{h}(k)\neq 0$ for all units of $\mathbb{Z}_q$, of which there
  are $\phi(q)$. In the latter case only if $\phi(q)\leq 2$ is it the
  case that $\wh{f}(k)\wh{h}(k)=0$ for $k\not\in\{1,-1\}$. The
  former case by Lemma~\ref{support} reduces to the latter with $q$ replaced by $d$.  
\end{pf}

The only choices for $q\ge 3$ and $S, S'\subseteq\mathbb{Z}_q$ are up to
exceptions trivial by Corollary~\ref{when}
given by the following
theorem, whose proof is a simple matter of substituting in the
expressions provided by Lemma~\ref{bias S} and Lemma~\ref{Bias Sigma
  Delta}.

\begin{thm}\label{general} Let $q\in\{3,4,6\}$. Suppose $A,B\,(,C)\subseteq E$ are chosen uniformly at
  random and $\Delta$ is the event that $A\bigtriangleup B$ (and
  $B\bigtriangleup C$) is eulerian. Suppose further that $S, S'\subseteq
    \mathbb{Z}_q$ and ${\rm supp}(\wh{1_S}-\wh{1_{S'}})=\{1,-1\}$
    (or possibly $\{1,-1, q/2\}$ for the third case of the following statement). 
If $\Sigma$ is the event
  $|A|\pm|B|\,(+|C|)\in S\,(\bmod\,q)$, $\Sigma'$ is the event
  $|A|\pm|B|\,(+|C|)\in S'\,(\bmod\,q)$ and ${\rm
  Bias}(\Sigma)\neq{\rm Bias}(\Sigma')$, then 
$$\frac{{\rm Bias}(\Sigma\,|\,\Delta)-{\rm Bias}(\Sigma'\,|\,\Delta)}{{\rm Bias}(\Sigma)-{\rm Bias}(\Sigma')}$$
$$=\begin{cases} 2^{r(G)}\left(1+\cos\frac{2\pi
      }{q}\right)^{-|E|}\left(\cos\frac{2\pi}{q}\right)^{r(G)}\left(\cos\frac{2\pi
      }{q}-1\right)^{n(G)}T\left(\!G;\frac{1}{\cos\frac{2\pi
        }{q}},\frac{1+\cos\frac{2\pi }{q}}{1-\cos\frac{2\pi}{q}}\right)\\
 2^{r(G)}\left(1+\cos\frac{2\pi
      }{q}\right)^{-|E|}\left(\cos\frac{2\pi}{q}-1\right)^{n(G)}T\left(\!G;\cos\frac{2\pi}{q},\frac{\cos\frac{2\pi}{q}+1}{\cos\frac{2\pi}{q}-1}\right) \\
 2^{r(G)}\left(1+\cos\frac{2\pi}{q}\right)^{-|E|}\left(\cos\frac{2\pi}{q}-1\right)^{n(G)}T\left(\!G;2\cos\frac{2\pi}{q}-1,\frac{\cos\frac{2\pi}{q}+1}{\cos\frac{2\pi}{q}-1}\right)\end{cases}$$
according as
 $$\Sigma=\begin{cases}\{A,B\subseteq E:|A|\!-\!|B|\in S\,(\bmod\,q)\}\\
\{A,B\subseteq E:|A|\!+\!|B|\in S\,(\bmod\,q)\}\\
\{A,B,C\subseteq E:|A|\!+\!|B|\!+\!|C|\in S\,(\bmod\,q)\}\end{cases},$$
$$\Sigma'=\begin{cases}\{A,B\subseteq E:|A|\!-\!|B|\in S'\,(\bmod\,q)\}\\
\{A,B\subseteq E: |A|\!+\!|B|\in S'\,(\bmod\,q)\}\\
\{A,B,C\subseteq E: |A|\!+\!|B|\!+\!|C|\in S'\,(\bmod\,q)\}\end{cases}.$$

\end{thm}
Taking $q$ even, $|S|=q/2$ and $S'=\mathbb{Z}_q\setminus S$, Theorem~\ref{general} gives ${\rm Bias}(\Sigma\,|\,\Delta)/{\rm Bias}(\Sigma)$ as a Tutte
polynomial evaluation.

So what choices of $S$ and $S'$ fulfil the conditions of Theorem~\ref{general}?
The answer is to be found in the theorems of sections~\ref{square root of unity} to \ref{sixth root of unity}, which are immediate
corollaries of Lemmas~\ref{bias S} and \ref{Bias Sigma Delta} and Theorem~\ref{general}. 

\subsubsection{Evaluations for $q=2$}\label{square root of unity}

\begin{prop} \label{squares with q=2}
Suppose that $A,B\subseteq E$ are subgraphs of $G$ chosen uniformly at
random. Then the event that $|A|+|B|$ is even is correlated with the
event $\Delta$ that $A\bigtriangleup B$ is eulerian as follows: 
$${\rm Bias}(\,|A|\!+\!|B|\equiv 0\,(\mbox{\rm mod}\, 2)\,\mid\,\Delta)=(-1)^{r(G)}T(G;-1,0)=2^{-k(G)}P(G;2).$$
\end{prop}
\begin{pf} In Lemma~\ref{Bias Sigma Delta} take $S=\{0\}$, $\Sigma=\{A,B\subseteq E:|A|+|B|\equiv
  0(\bmod 2)\}$ and $\Delta$ the event that $A\bigtriangleup B$ is
  eulerian. Then, by the result of that lemma, ${\rm
    Bias}(\Sigma\,|\,\Delta)=\wh{f}(0)+\wh{f}(1)-1$, where
  $\wh{f}(0)=1$ and $\wh{f}(1)=2^{-n(G)}(-1)^{r(G)}2^{r(G)}T(G;-1,0)=2^{-k(G)}P(G;2)$.   
\end{pf} 

 Of course the correlation between parity and eulerian symmetric difference in
Proposition~\ref{squares with q=2} can be seen immediately by considering the identity 
$$|A\bigtriangleup B|+2|A\cap B|=|A|+|B|.$$
 Eulerian subgraphs are all of even size if and
only if $G$ is bipartite (no odd cycles). 
 Given the event $\Delta$ that $A\bigtriangleup B$
is eulerian the parity of $|A|+|B|$ must be even when $G$ is
bipartite.
Otherwise, if $G$ is not bipartite half the eulerian
subgraphs are even, half odd, and so the parity of $|A|+|B|$ is
equally likely to be even or odd given $\Delta$.

For three subgraphs $A,B,C\subseteq
E$ of $G$, if $A\bigtriangleup  B, B\bigtriangleup  C$ are eulerian then so is $C\bigtriangleup A$.
From the identity
$$|A\bigtriangleup B|+|B\bigtriangleup C|+|C\bigtriangleup A|=2(|A|+|B|+|C|)-2(|A\cap
B|+|B\cap C|+|C\cap A|),$$
it seems difficult to tell whether there might be any correlation between the
event that $A\bigtriangleup  B, B\bigtriangleup  C$ are eulerian and
some condition on $|A|+|B|+|C|$.

\begin{thm}\label{three uncorr even}
Suppose $A,B,C\subseteq E$ are subgraphs of $G$ chosen uniformly at
random. Then the event
that $|A|+|B|+|C|$ is even is uncorrelated with the event $\Delta$ that
$A\bigtriangleup  B$, $B\bigtriangleup  C$ are eulerian, i.e., $${\rm
  Bias}(\,|A|\!+\!|B|\!+\!|C|\equiv 0\!\!\pmod 2\,|\,\Delta)=0.$$
\end{thm}
\begin{pf} Take $\Sigma=\{A,B,C\subseteq E:|A|+|B|+|C|\equiv 0(\bmod\, 2)\}$ and $\Delta$ the event that $A\bigtriangleup B$ and
  $B\bigtriangleup C$ are both
  eulerian. By Lemma~\ref{Bias Sigma Delta}, ${\rm
    Bias}(\Sigma\,|\,\Delta)=\wh{f}(0)+\wh{f}(1)-1$, where
  $\wh{f}(0)=1$ and $\wh{f}(1)=0$.   
\end{pf} 

However, we shall see that the residue of $|A|+|B|+|C|$ modulo $3, 4$ and
$6$ does have a
bearing on the event that $A\bigtriangleup B, B\bigtriangleup C$ are
eulerian. 

\subsubsection{Evaluations for $q=3$}\label{cube root of unity}

The evaluations of the Tutte polynomial obtained for $q=3$ are, unlike
the cases $q=2, 4$ and $6$, at
points without other more familiar combinatorial interpretations.

\begin{thm}\label{A - B q=3}
 Let $A,B\subseteq E$ be chosen uniformly at random and let $\Delta$
 be the event that $A\bigtriangleup B$ is eulerian. Then
$${\rm Bias}(\,|A|\equiv |B|\!+\!1\,(\bmod\,
3)\,\mid\,\Delta)={\rm Bias}(\,|A|\equiv|B|\!+\!2\,(\bmod\, 3)\,\mid\,\Delta)$$
and 
$$\frac{{\rm Bias}(\,|A|\equiv |B|\,(\bmod\,
  3)\,\mid\,\Delta)-{\rm Bias}(\,|A|\equiv |B|\!+\!1\,(\bmod\, 3)\,\mid\,\Delta)}{{\rm Bias}(\,|A|\equiv
  |B|\,(\bmod\, 3)\,)-{\rm Bias}(\,|A|\equiv
  |B|\!+\!1\,(\bmod\, 3)\,)}$$
$$=(-2)^{r(G)}3^{n(G)}T(G;-2,\frac{1}{3}).$$
\end{thm}
\begin{pf}
Lemma~\ref{Bias Sigma Delta} with $S=\{1\}$ yields 
$${\rm Bias}(|A|-|B|\equiv 1(\bmod\, 3)\,|\,\Delta)=\frac23(1+\wh{f}(1)e^{2\pi
  i/3}+\wh{f}(2)e^{4\pi i/3})-1,$$
where $f(\ell)=\mathbb{P}(|A|-|B|\equiv\ell (\bmod\,3)\,|\,\Delta)$. (We  define a
character $\chi$ on $\mathbb{Z}_3$ by setting $\chi(1)=e^{2\pi i/3}$, and the Fourier transform is
defined by $\wh{f}(k)=f(0)+e^{4k\pi i/3}f(1)+e^{2k\pi i/3}f(2)$. Thus
$\wh{1_1}=1_0+e^{4\pi i/3}1_1+e^{2\pi i/3}1_2$ and
$\wh{1}_2=\overline{\wh{1}_1}$.)

 With $S=\{2\}$ the same lemma yields 
$${\rm Bias}(|A|-|B|\equiv 2(\bmod 3)\,|\,\Delta)  = \frac23(1+\wh{f}(1)e^{4\pi
  i/3}+\wh{f}(2)e^{2\pi i/3})-1.$$
Since Lemma \ref{Bias Sigma Delta} also tells us that  $\wh{f}(1)=\wh{f}(2)$, the first statement of the theorem
is established.

By Theorem~\ref{general} with $S=\{0\}, S'=\{1\}$ (for which
$\wh{1_S}-\wh{1_{S'}}=(1-e^{4\pi i/3})1_1+(1-e^{2\pi i/3})1_2$),
$\Sigma=\{A,B\subseteq E:|A|-|B|\equiv 0(\bmod 3)\}$ and
$\Sigma'=\{A,B\subseteq E: |A|-|B|\equiv 1 (\bmod 3)\}$,
$$\frac{{\rm Bias}(\,|A|-|B|\equiv 0\,(\mbox{\rm mod}\,
  3)\,\mid\,\Delta)-{\rm Bias}(\,|A|-|B|\equiv 1\,(\mbox{\rm
    mod}\, 3)\,\mid\,\Delta)}{{\rm Bias}(\,|A|-|B|\equiv 0\,(\mbox{\rm mod}\, 3)\,)-{\rm Bias}(\,|A|-|B|\equiv
  1\,(\mbox{\rm mod}\, 3)\,)}$$
$$=2^{r(G)}(\frac12)^{-|E|}(-\frac12)^{r(G)}(-\frac32)^{n(G)}T(G;-2,\frac13).$$ 
\end{pf}

\begin{thm}\label{cube root A B}
 Let $A,B\subseteq E$ be chosen uniformly at random and let $\Delta$
 be the event that $A\bigtriangleup B$ is eulerian. Then
$${\rm Bias}(\,|A|\!+\!|B|\equiv |E|\!+\!1\,(\mbox{\rm mod}\,
3)\,\mid\,\Delta)={\rm Bias}(\,|A|\!+\!|B|\equiv |E|\!+\!2\,(\mbox{\rm mod}\, 3)\,\mid\,\Delta)$$
and 
$$\frac{{\rm Bias}(\,|A|\!+\!|B|\equiv |E|\,(\mbox{\rm mod}\,
  3)\,\mid\,\Delta)-{\rm Bias}(\,|A|\!+\!|B|\equiv |E|\!+\!1\,(\mbox{\rm
    mod}\, 3)\,\mid\,\Delta)}{{\rm Bias}(\,|A|\!+\!|B|\equiv
  |E|\,(\mbox{\rm mod}\, 3)\,)-{\rm Bias}(\,|A|\!+\!|B|\equiv
  |E|\!+\!1\,(\mbox{\rm mod}\, 3)\,)}$$
$$=4^{r(G)}(-3)^{n(G)}T(G;-\frac{1}{2},-\frac{1}{3}).$$\end{thm}

\begin{pf}  Lemma~\ref{Bias Sigma Delta} with $S=\{|E|\!+\!1\}$ yields 
$${\rm Bias}(|A|+|B|\equiv |E|\!+\!1(\bmod 3)\,|\,\Delta)=\frac23(1+\wh{f}(1)e^{2\pi
  i(|E|+1)/3}+\wh{f}(2)e^{4\pi i(|E|+1)/3})-1,$$
where $f(\ell)=\mathbb{P}(|A|+|B|\equiv\ell (\bmod
3)\,|\,\Delta)$. With $S=\{|E|\!+\!2\}$ the same lemma yields 
$${\rm Bias}(|A|+|B|\equiv |E|\!+\!2(\bmod 3)\,|\,\Delta)  = \frac23(1+\wh{f}(1)e^{2\pi i(|E|+2)/3}+\wh{f}(2)e^{4\pi i(|E|+2)/3})-1.$$
Lemma~\ref{Bias Sigma Delta} tells us that $\wh{f}(1)=e^{2\pi
  i|E|/3}\wh{f}(2)$, and the first statement of the theorem follows
with both biases equal to $\wh{f}(1)(e^{2\pi i(|E|+1)/3}+e^{2\pi i(|E|+2)/3})-\frac13$.

The second statement of the theorem follows from Theorem~\ref{general}
upon taking $S=\{|E|\}, S'=\{|E|\!+\!1\}$,
$\Sigma=\{A,B\subseteq E:|A|\!+\!|B|\equiv |E|(\bmod{3})\,\}$ and
$\Sigma'=\{A,B\subseteq E: |A|\!+\!|B|\equiv |E|\!+\!1(\bmod{3})\,\}$.
\end{pf}

Note that by equation~(\ref{eul reciprocal}) at the beginning of this
section, if $G$ is eulerian then the evaluations of the Tutte polynomial in Theorem~\ref{A -
  B q=3} and Theorem~\ref{cube
  root A B} are equal. Indeed,  $|E\setminus A|+|B|\equiv\pm
  1\,(\mbox{\rm mod}\, 3)$ if and only if $|A|-|B|\equiv|E|\mp
  1\,(\mbox{\rm mod}\, 3)$, and if $G$ is eulerian
  then a subgraph $A$ is eulerian if and only if its complement
  $E\setminus A$ is eulerian. 
  
\begin{thm}\label{cube root A B C} Let $A,B,C\subseteq E$ be chosen uniformly at random and let $\Delta$ be the event that $A\bigtriangleup B, B\bigtriangleup
  C$ are both eulerian. Then 
$${\rm Bias}(\,|A|\!+\!|B|\!+\!|C|\equiv 1\,(\bmod 3)\,\mid\,\Delta)={\rm
  Bias}(\,|A|\!+\!|B|\!+\!|C|\equiv 2\,(\bmod 3)\,\mid\,\Delta),$$ and
$$\frac{{\rm Bias}(\,|A|\!+\!|B|\!+\!|C|\equiv 0\,(\bmod 3)\,\mid\,\Delta)-{\rm Bias}(\,|A|\!+\!|B|\!+\!|C|\equiv 1\,(\bmod 3)\,\mid\,\Delta)}{{\rm
    Bias}(\,|A|\!+\!|B|\!+\!|C|\equiv 0\,(\bmod 3)\,)-{\rm Bias}(\,|A|\!+\!|B|\!+\!|C|\equiv 1\,(\bmod 3)\,)}$$
$$=4^{r(G)}(-3)^{n(G)}T(G;-2,-\frac{1}{3}).$$
\end{thm}

\begin{pf}
 Lemma~\ref{Bias Sigma Delta} with $S=\{1\}$ yields
$${\rm Bias}(|A|\!+\!|B|\!+\!|C|\equiv
 1(\bmod{3})\,|\,\Delta)=\frac23(\wh{f}(1)e^{2\pi
   i/3}+\wh{f}(2)e^{4\pi i/3})-\frac13$$
where $f(\ell)=\mathbb{P}(|A|\!+\!|B|\!+\!|C|\equiv\ell (\bmod
3)\,|\,\Delta)$.   Lemma~\ref{Bias Sigma Delta} with $S=\{2\}$ yields
$${\rm Bias}(|A|\!+\!|B|\!+\!|C|\equiv
 2(\bmod{3})\,|\,\Delta)=\frac23(\wh{f}(1)e^{4\pi
   i/3}+\wh{f}(2)e^{2\pi i/3})-\frac13.$$
From  Lemma~\ref{Bias Sigma Delta} it is also found that
$\wh{f}(1)=\wh{f}(2)$ and the first statement of the theorem follows.

Clearly then ${\rm Bias}(\,|A|\!+\!|B|\!+\!|C|\equiv 0(\bmod\,3)\,)\neq
{\rm Bias}(\,|A|\!+\!|B|\!+\!|C|\equiv 1(\bmod\,3)\,)$ and the second
statement of the theorem results from Theorem~\ref{general} upon taking $S=\{0\},$  $S'=\{1\}$,
  $\Sigma=\{A,B,C\subseteq E:|A|\!+\!|B|\!+\!|C|\equiv
 0(\bmod\,3)\}$ and $\Sigma'=\{A,B,C\subseteq E:|A|\!+\!|B|\!+\!|C|\equiv
 1(\bmod\,3)\}$.
\end{pf}

\subsubsection{Evaluations for $q=4$}\label{fourth root of unity}

\begin{thm}\label{A - B with q=4}
Choosing $A,B\subseteq E$ uniformly at random, the event that $|A|-|B|\equiv
0\,\mbox{\rm or}\, 1\,(\bmod\, 4)$ (i.e.,
$\left\lfloor\frac{|A|-|B|}{2}\right\rfloor$ is even) is correlated with the event $\Delta$ that $A\bigtriangleup 
B$ is eulerian:
$$\frac{{\rm Bias}(\,|A|\!-\!|B|\equiv
  0,1\,(\bmod 4)\,\mid\,\Delta)}{{\rm
    Bias}(|A|\!-\!|B|\equiv 0,1\,(\bmod\, 4)\,)}=2^{r(G)}.$$ 
\end{thm}
\begin{pf}
In Lemma~\ref{bias S} take $S=\{0,1\}$, for which
$\wh{1_S}(k)=1+i^{-k}$, and calculate 
$${\rm Bias}(|A|\!-\!|B|\equiv 0,1(\bmod{4}))=2^{-1-|E|}\big[2^{|E|}\cdot
2+1\cdot(1+i)+0\cdot 0+1\cdot (1-i)\big]-1=2^{-|E|}.$$

That this is non-zero allows us to apply Theorem~\ref{general}, in which we take $S=\{0,1\}$,
$S'=\{2,3\}=\mathbb{Z}_4\setminus S$, $\Sigma=\{A,B\subseteq
E:|A|\!-\!|B|\equiv 0,1(\bmod\,4)\}$ and $\Sigma'=\{A,B\subseteq
E:|A|\!-\!|B|\equiv 2,3(\bmod\,4)\}$. This yields (recalling from the footnote on page~\pageref{fn} how to
deal with division by zero in Tutte polynomial evaluations)
$$\frac{{\rm Bias}(\,|A|\!-\!|B|\equiv
  0,1(\bmod\,4)\,\mid\,\Delta)-{\rm Bias}(|A|\!-\!|B|\equiv 2,3(\bmod{4})\,|\,\Delta)}{{\rm
    Bias}(|A|\!-\!|B|\equiv 0,1(\bmod\,4)\,)-{\rm
    Bias}(|A|\!-\!|B|\equiv 2,3(\bmod\,4)\,)}=2^{r(G)}.$$
\end{pf}

\begin{thm} \label{squares with q=4}
Choosing $A,B\subseteq E$ uniformly at random, the event that $|A|+|B|\equiv
0\,\mbox{\rm or}\, 1\,(\mbox{\rm mod}\, 4)$ (i.e.,
$\left\lfloor\frac{|A|+|B|}{2}\right\rfloor$ is even) is correlated with the event $\Delta$ that $A\bigtriangleup 
B$ is eulerian in the following way: 
$$\frac{{\rm Bias}(\,|A|\!+\!|B|\equiv
  0,1(\bmod\,4)\,\mid\,\Delta)}{{\rm
    Bias}(|A|\!+\!|B|\equiv 0,1(\bmod\, 4))}=2^{r(G)}F(G;2).$$ 
\end{thm}

\begin{pf}
In Lemma~\ref{bias S} take $S=\{0,1\}$, for which
$\wh{1_S}(k)=1+i^{-k}$, and calculate 
\begin{align*}
{\rm Bias}(|A|+\!|B|\equiv 0,1(\bmod\,4)) & =
2^{-1-|E|}\big[2^{|E|}\cdot 2+i^{-|E|}\cdot(1+i)+0\cdot 0+i^{|E|}\cdot (1-i)\big]-1\\
& =(-1)^{\lfloor|E|/2\rfloor}2^{-|E|}.
\end{align*} 

We can now apply Theorem~\ref{general}, taking $S=\{0,1\}$,
$S'=\{2,3\}=\mathbb{Z}_4\setminus S$, $\Sigma=\{A,B\subseteq
E:|A|\!+\!|B|\equiv 0,1(\bmod\,4)\}$ and $\Sigma'=\{A,B\subseteq
E:|A|\!+\!|B|\equiv 2,3(\bmod\,4)\}$. The formula in Theorem~\ref{general} yields
$$\frac{{\rm Bias}(\,|A|\!+\!|B|\equiv
  0,1(\bmod{4})\,\mid\,\Delta)-{\rm Bias}(|A|\!+\!|B|\equiv 2,3(\bmod\,4)\,|\,\Delta)}{{\rm
    Bias}(|A|\!+\!|B|\equiv 0,1(\bmod\,4))-{\rm
    Bias}(|A|\!+\!|B|\equiv
  2,3(\bmod\,4))}$$
$$=2^{r(G)}(-1)^{n(G)}T(G;0,-1).$$
\end{pf}

Theorem~\ref{squares with q=4} says that if $G$ is not eulerian then the event that $A\bigtriangleup 
B$ is eulerian removes from the
 parity of $\lfloor\frac{|A|+|B|}{2}\rfloor$ its original bias of
$(-1)^{\lfloor|E|/2\rfloor}2^{-|E|}$ towards being even. Otherwise, when $G$ is eulerian, this bias is accentuated by a
factor of $2^{r(G)}$.
Theorem~\ref{squares with q=4} has a counterpart in Theorem~\ref{sum cubes with q=6} in section~\ref{sixth root of unity} below.

\begin{thm} \label{A + B + C with q=4}
Suppose $|E|\not\equiv 1\pmod{4}$. Then, for $A,B,C\subseteq E$ chosen
uniformly at random, the event that $|A|\!+\!|B|\!+\!|C|\equiv
0\,\mbox{\rm or}\, 1\,\pmod 4$ (i.e.,
$\left\lfloor\frac{|A|+|B|+|C|}{2}\right\rfloor$ is even) is correlated with the event $\Delta$ that $A\bigtriangleup 
B$ and $B\bigtriangleup C$ are eulerian as follows: 
$$\frac{{\rm Bias}(\,|A|\!+\!|B|\!+\!|C|\equiv
  0,1\,(\bmod\, 4)\,\mid\,\Delta)}{{\rm
    Bias}(|A|\!+\!|B|\!+\!|C|\equiv 0,1\,(\bmod\, 4))}=2^{r(G)}(-1)^{n(G)}T(G;-1,-1).$$ 
If $|E|\equiv 1\pmod 4$ then 
$${\rm  Bias}(|A|\!+\!|B|\!+\!|C|\equiv
0,1(\bmod\, 4)\,)=0={\rm Bias}(\,|A|\!+\!|B|\!+\!|C|\equiv
  0,1\,(\bmod\, 4)\,\mid\,\Delta).$$
\end{thm}
\begin{pf}
In Lemma~\ref{bias S} take $q=4$ and $S=\{0,1\}$, for which
$\wh{1_S}(k)=1+i^{-k}$, and calculate 
\begin{align*}
{\rm Bias}(|A|\!+\!|B|\!+\!|C|\equiv 0,1(\bmod{4})) & = 2^{-1-\frac32|E|}\big[2^{\frac32|E|}\cdot
2+i^{-\frac32|E|}\cdot(1\!+\!i)+i^{\frac32|E|}\cdot
(1\!-\!i)\big]-1\\
& =2^{-\frac32|E|}{\rm Re}\big[i^{\frac32|E|}(1-i)\big]\\
& ={\begin{cases} 2^{-3|E|/2} & |E|\equiv 0,6\pmod 8\\
0 & |E|\equiv 1,5\pmod 8\\
-2^{-3|E|/2} & |E|\equiv 2,4\pmod 8\\
2^{(1-3|E|)/2} & |E|\equiv 3\pmod 8\\
-2^{(1-3|E|)/2} & |E|\equiv 7\pmod 8\end{cases}.}\end{align*}
When $|E|\not\equiv 1 (\bmod\, 4)$ we can apply Theorem~\ref{general} and
the result follows by a straightforward calculation.

When ${\rm Bias}(\,|A|\!+\!|B|\!+\!|C|\equiv 0,1 (\bmod\,4) )=0$ we use
Lemma~\ref{Bias Sigma Delta} with $q=4, S=\{0,1\}$, to calculate that
\begin{align*}
{\rm Bias}(|A|\!+\!|B|\!+\!|C|\equiv 0,1(\bmod\,4)\,|\,\Delta) & =
2{\rm Re}\big[2^{-n(G)-\frac12|E|}i^{\frac32|E|}(-1)^{n(G)}T(G;-1,-1)(1-i)\big]\\
& =0\hspace{1cm}\mbox{\rm when }|E|\equiv 1 (\bmod\, 4).\end{align*} 

\end{pf}

From~\cite{RR78},
$2^{r(G)}(-1)^{n(G)}T(G;-1,-1)=(-2)^{r(G)+\dim(\mathcal{C}_2\cap\mathcal{C}_2^\perp)}$,
where $\mathcal{C}_2\cap\mathcal{C}_2^\perp$ is the bicycle space of
$G$, comprising subgraphs which are both eulerian and bipartite.


\subsubsection{Evaluations for $q=6$}\label{sixth root of unity}
 
When $q=6$ and $S=\{0,1,2\}$ the expression ${\rm
  Bias}(\Sigma\,\mid\,\Delta)/{\rm Bias}(\Sigma)$ is only equal to an evaluation of the Tutte
polynomial at a single point when $\Sigma=\{A,B,C\subseteq E:|A|\!+\!|B|\!+\!|C|\equiv 0,1,2\,(\bmod\, q)\}$. This is due to the formula for ${\rm
  Bias}(\Sigma\,\mid\,\Delta)$ given by Lemma~\ref{Bias Sigma Delta}
  and the fact that ${\rm supp}(\wh{1_S})=\{0,\pm 1, 3\}$.
  For example, when $\Sigma$ is the event
$\{A,B\subseteq E:|A|\!+\!|B|\equiv 0,1,2\,(\bmod\, 6)\,\}$,
  evaluations of the Tutte polynomial at the two points
  $(\frac{1}{2},-3)$ and $(0,-1)$ would be involved.

\begin{thm}\label{sum cubes with q=6} Suppose that $A,B,C\subseteq E$ are
  chosen uniformly at random and $\Delta$ is the event that $A,B,C$ have
  pairwise eulerian differences. Then the event that
  $|A|+|B|+|C|\equiv 0,1,2\,(\bmod\, 6)$ (i.e., $\left\lfloor\frac{|A|+|B|+|C|}{3}\right\rfloor$ is
even) is correlated with $\Delta$ as follows:
$$\frac{{\rm Bias}(\,|A|\!+\!|B|\!+\!|C|\equiv 0,1,2\,(\bmod\, 6)\,\mid\,\Delta)}{{\rm Bias}(\,|A|\!+\!|B|\!+\!|C|\equiv 0,1,2\,(\bmod\, 6)\,)}=3^{-|E|}4^{r(G)}F(G;4).$$
\end{thm}
\begin{pf}
In Lemma~\ref{bias S} take $q=6$ and $S=\{0,1,2\}$, for which we calculate that
$\wh{1_S}=31_0+1_3-2e^{2\pi i/3}1_1-2e^{-2\pi i/3}1_5$. By
Lemma~\ref{bias S}, 
\begin{align*}
{\rm Bias}(\,|A|\!+\!|B|\!+\!|C|\equiv 0,1,2\,(\bmod\,6)\,) & =\frac{2}{6}\big[\wh{g}(1)\ol{\wh{1_S}(1)}+\wh{g}(5)\ol{\wh{1_S}(5)}\big]\\
& = \frac23{\rm Re}\big[\wh{g}(1)\ol{\wh{1_S}(1)}\big]
\end{align*}
where $\wh{g}(1)=2^{-\frac32|E|}e^{-2\pi
  i\frac32|E|/6}(1+\cos\frac{2\pi}{6})^{|E|}=2^{-\frac32|E|}e^{-|E|\pi i /2}(1+\frac12)^{\frac32|E|}$. 
Hence 
\begin{align*}
{\rm Bias}(\,|A|\!+\!|B|\!+\!|C|\equiv 0,1,2\,(\bmod\,6)\,) & =
\frac23{\rm Re}\big[e^{-\pi i|E|/2}\big(\frac34\big)^{\frac32|E|}(2e^{-2\pi i/3})\big]\\
 & =\big(\frac34\big)^{\frac32|E|-1}{\rm Re}\big[i^{-|E|}e^{-2\pi
    i/3}\big]\neq 0
\end{align*}
and we can apply Theorem~\ref{general} from which the result follows
by routine calculation.
\end{pf}

\section{A parity criterion for proper vertex colourings}\label{poset}

In the final three sections of this article we
    need some further identities from finite Fourier analysis for complete weight enumerators
    (which include Hamming weight enumerators as specialisations). 

Let $Q$ be a commutative ring with a generating character such as $\mathbb{Z}_q$ or
    $\mathbb{F}_q$, and let $f$ be a function
    $f:Q\rightarrow\mathbb{C}$. 

The {\em complete weight enumerator} of 
a subset $\mathcal{S}$ of $Q^E$ is defined by
$${\rm cwe}(\mathcal{S};f)=\sum_{x\in\mathcal{S}}\prod_{e\in
  E}f(x_e).$$
When $f=t1_0+1_{Q\setminus 0}$ the complete weight enumerator is
the Hamming weight enumerator ${\rm hwe}(\mathcal{S};t)$.
The MacWilliams duality theorem for complete weight enumerators is a
consequence of the Poisson summation formula and states
that when $\mathcal{S}$ is a $Q$-submodule of $Q^E$
\be\label{McW cwe}{\rm cwe}(\mathcal{S};f)=\frac{1}{|\mathcal{S}^\perp|}{\rm
  cwe}(\mathcal{S}^\perp;\wh{f}).\ee

The following generalises the first two identities of Lemma~\ref{sum cubes} and is proved for example in \cite{CCC07}.
 
\begin{lem}\label{cwe sum} 
Let $Q^E$ be a commutative ring with a generating
   character. For $Q$-submodule $\mathcal{S}$ of $Q^E$ and functions
   $f,g:Q\rightarrow\mathbb{C}$,
   $$\sum_{\mathcal{S}+z\in Q^E/\mathcal{S}}{\rm
  cwe}(\mathcal{S}+z;f)\ol{{\rm
   cwe}(\mathcal{S}+z;g)}=\frac{1}{|\mathcal{S}^\perp|}{\rm
  cwe}(\mathcal{S}^\perp;\wh{f}\cdot\ol{\wh{g}}).$$
\end{lem}

Let $\mathcal{C}$ be the set of
$Q$-flows of $G$ and its orthogonal $\mathcal{C}^\perp$ the set of
$Q$-tensions of $G$. A partial order $\leq$ on $Q^E$ is defined by
$x\leq y$ if and only if $x_e\in\{0,y_e\}$ for all $e\in E$. (For
$Q=\mathbb{F}_2$ the order $\leq$ is set inclusion.) This makes the
poset on $Q^E$ the direct product of the poset $P$ on $Q$ defined by setting
$0$ below all the non-zero elements of $Q$ and all pairs of non-zero
elements incomparable. Thus the M\"{o}bius function of the poset $P^E=(Q^E,\leq)$ is
defined by $\mu(x,y)=(-1)^{|y|-|x|}$. (See for example \cite{St97} for
background on posets.) For a function
$f:Q^E\rightarrow\mathbb{C}$, define $\mu f:Q^E\rightarrow\mathbb{C}$ by
$$\mu f(y)=\sum_{x\leq y}\mu(x,y)f(x)=\sum_{x\leq
  y}(-1)^{|y|-|x|}f(x).$$

 \begin{lem}\label{mu1C not zero}
Let $Q$ be a ring with a generating character $\chi$.
If $\mathcal{C}$ is a $Q$-submodule of $Q^E$ and
$\mathcal{C}^\perp$ 
 its orthogonal space then
$$\mu 1_\mathcal{C}(y)=\frac{1}{|\mathcal{C}^\perp|}\sum_{x\in\mathcal{C}^\perp}\prod_{e\in E}(\chi(x_ey_e)-1).$$
\end{lem}
\begin{pf}
\begin{align*}\mu 1_{\mathcal{C}}(y) &
  =\sum_{x\in\mathcal{C}}\prod_{e\in E}(1_{y_e}-1_0)(x_e)\\
 & = \frac{1}{|\mathcal{C}^\perp|}\sum_{x\in\mathcal{C}^\perp}(\ol{\chi}_{y_e}-1)(x_e),\end{align*}
the latter equality by identity \eqref{McW cwe}, and since the
left-hand side is real $\ol{\chi}_{y_e}(x_e)$ can be replaced by its
conjugate $\chi_{y_e}(x_e)=\chi(x_ey_e)$.
\end{pf}
\begin{lem} \label{sum nz C perp}
$$\sum_{y: \forall_{e\in E}\, y_e\neq 0}\mu
  1_\mathcal{C}(y)=(-1)^{|E|}|\mathcal{C}|{\rm
   hwe}(\mathcal{C}^\perp;0).$$\end{lem}

\begin{pf}
\begin{align*}\sum_{y:\, \forall_{e\in E}\, y_e\neq 0}\mu
  1_\mathcal{C}(y) & = \sum_{y:\, \forall_{e\in E}\, y_e\neq 0}\sum_{x\leq y,\, x\in\mathcal{C}}(-1)^{|y|-|x|}\\
 &
 =\sum_{x\in\mathcal{C}}(q-1)^{|E|-|x|}(-1)^{|E|-|x|},\end{align*}
reversing the order of summation and using $\#\{y\in(Q\setminus
0)^E:x\leq y\}=(q-1)^{|E|-|x|}$, whence
\begin{align*}\sum_{y:\, \forall_{e\in E}\, y_e\neq 0}\mu
  1_\mathcal{C}(y) & =\sum_{x\in\mathcal{C}}(1-q)^{|E|-|x|}\\
 & ={\rm hwe}(\mathcal{C};1-q)=\frac{(-q)^{|E|}}{|\mathcal{C}^\perp|}{\rm
   hwe}(\mathcal{C}^\perp;0),\end{align*}
and finally $q^{|E|}/|\mathcal{C}^\perp|=|\mathcal{C}|$.
\end{pf}

The following is a variation on, and mild generalisation of, Theorem~1.2 in \cite{Onn04}.

\begin{cor}\label{cond Qtensions}
Suppose $G$ is a graph and $Q$ is a ring of order $q$ with a
generating character. Let $\mathcal{C}$ be the set
of $Q$-flows of $G$ and
$\mathcal{C}^\perp$ the set of $Q$-tensions of $G$. Then $P(G;q)\neq
0$ if and only if there exists 
$y\in (Q\setminus 0)^E$ such that $\mu 1_\mathcal{C}(y)\neq 0$,
i.e., such that
$$\sum_{x\leq y, \,
  x\in\mathcal{C}}(-1)^{|x|}\neq 0.$$
\end{cor}
\begin{pf} From Lemma~\ref{mu1C not zero}, if $\chi(x_ey_e)\neq 1$ for all $e\in E$ then
  $x_e\neq 0$ for all $e\in E$. The converse follows from Lemma~\ref{sum nz C perp}. 
\end{pf}
A dual to Corollary~\ref{cond Qtensions} giving a criterion for
$F(G;q)\neq 0$ results by taking $\mathcal{C}$ to be the set of
$Q$-tensions. Corollary~\ref{cond Qtensions} was proved for $q=3$ and
generalised in a different direction by Alon and Tarsi~\cite{AT92} by
considering $\mathbb{Z}_q$-flows taking values in $\{0,\pm 1\}$ only. The latter
for $q$ greater than the maximum degree of $G$ are in bijective
correspondence with partial eulerian orientations
of $G$ (and for $q=3$ the same is true for $4$-regular graphs). See
also \cite{Tarsi} and
section~\ref{4 regular} below. 

From Corollary~\ref{cond Qtensions} comes the familiar fact that $P(G;2)\neq 0$ if and only if the number of eulerian
subgraphs of $G$ of even size differs from those with odd size. More
interestingly, $P(G;4)\neq 0$ if and only if there is
$y\in\mathbb{F}_4^E$ such that $\mu 1_\mathcal{C}(y)\neq 0$, i.e., the difference between the number of
$\mathbb{F}_4$-flows $\leq y$ of even support size and those $\leq y$
of odd support size is non-zero, where in this case
 \begin{align*}\mu
1_\mathcal{C}(y) & =4^{-r(G)}(-2)^{|E|}\#\{x\in\mathcal{C}^\perp:\forall_{e\in
  E}\; x_e\not\in\{0, y_e\}\}\\
 & = (-2)^{|E|-2|V|}\#\{z\in \mathbb{F}_4^V:\forall_{uv\in E}\; z_u+z_v\not\in
 \{0,y_e\}\}.\end{align*}
It may be that $\mu 1_\mathcal{C}(y)=0$ for some $y\in
(\mathbb{F}_4^\times)^E$ even though $P(G;4)\neq 0$, since it may be
impossible to avoid hitting the value $y_e$ for some edge $e$ in any
nowhere-zero $\mathbb{F}_4$-tension $x$ of $G$. (For example, the
triangle $K_3$ and $y_e=1$ for each edge $e$).

Similarly $P(G;4)\neq 0$ if and only if for some
$y\in(\mathbb{Z}_4\setminus 0)^E$ there is a disparity between the
number of $\mathbb{Z}_4$-flows $\leq y$ of even support size and those
$\mathbb{Z}_4$-flows $\leq y$ of odd support size, and here  
$$\mu 1_\mathcal{C}(y)=4^{-r(G)}\sum_{x\in\mathcal{C}^\perp}\prod_{e\in
  E}(i^{x_ey_e}-1).$$ 
If $y_e=2$ then in order for $x\in\mathcal{C}^\perp$ to contribute
a non-zero term it is necessary that $x_e\not\in\{0, 2\}$. This too may not be
possible for some $y$ (consider $K_3$ again with $y_e=2$ for each edge).
 

A translation of  Corollary~\ref{cond Qtensions} for $Q=\mathbb{F}_4$
into the language of correlations between events involving parity and
eulerian subgraphs runs as follows.

\begin{thm}\label{tripart} Suppose $X,Y,Z\subseteq E$ partition the edges of
  $G$ into three sets (not all of which need be non-empty). 
Choosing $A\subseteq X,
B\subseteq Y$ and $C\subseteq Z$ uniformly at random, let $\Sigma$ be the event that
$|A|+|B|+|C|\equiv 0\,(\mbox{\rm mod}\, 2)$ and $\Gamma$ the
event that $A\cup C$ and $C\cup B$ are both
eulerian.

Then ${\rm Bias}(\Sigma\,\mid\, \Gamma)\neq 0$ for some tripartition
$\{X,Y,Z\}$ of $E$ if and only if $P(G;4)\neq
0$. 
\end{thm}
Note that in contrast to the event $\Delta$ of section
\ref{H2 H4}, it does not follow that if $A\cup C,
C\cup B$ are eulerian then $A\cup
B$ is eulerian. Also, note that ${\rm Bias}(\Sigma)=0$ for any choice of $X,Y,Z$.
\begin{pf} 
We use Corollary~\ref{cond Qtensions} to show the auxilary result that $P(G;4)\neq 0$ if and only if there exist $X,Y\subseteq E$ with $X\cup Y=E$ and
\be\label{aux}\mathop{\sum_{\mbox{\rm \tiny eulerian } A\subseteq X,\, B\subseteq
      Y}}_{A\bigtriangleup B\subseteq X\bigtriangleup Y}(-1)^{|A\cup
    B|}\neq 0.\ee
We then take $A\setminus B$, $B\setminus A$ in \eqref{aux} for the $A$ and $B$ of the
theorem, $X\setminus Y, Y\setminus X$ in \eqref{aux} for the $X$ and $Y$ of the
theorem, and finally set $C=A\cap B$ and $Z=X\cap Y$. This is enough to
prove the theorem as stated, for \eqref{aux}, with $|(A\cup C)\cup (C\cup
B)|=|A|+|B|+|C|$, is now the assertion that
$$\mathbb{P}(\Sigma\cap\Gamma)-\mathbb{P}(\ol{\Sigma}\cap\Gamma)=2^{-3|E|}\sum_{\mbox{\rm
    \tiny eulerian }\, A\subseteq X, B\subseteq Y, C\subseteq
  Z}(-1)^{|A|+|B|+|C|}\neq 0.$$ 

Let $x,y\in\mathbb{F}_2^E$ be the indicator vectors of $X, Y$ and
 $z=(x,y)\in\mathbb{F}_2^E\times\mathbb{F}_2^E\cong\mathbb{F}_4^E$. 
 Define a partial
 order $\leq$ on
 $\mathbb{F}_4^E$ by setting $d\leq z$ if and only if $d_e\in\{0,
 z_e\}$. Then $d=(a,b)\leq z=(x,y)$ if and only if $A\subseteq X,
 B\subseteq Y$ and $A\bigtriangleup B\subseteq X\bigtriangleup Y$.
Note that $z_e\neq 0$ for all $e\in E$ if and only if $X\cup Y=E$. Denote
 by $|d|$ the Hamming weight of $d\in\mathbb{F}_4^E\cong\mathbb{F}_2^E\times\mathbb{F}_2^E$. If $d=(a,b)$ for
 $a,b\in\mathbb{F}_2^E$ the indicator vectors of $A,B\subseteq E$ then
 $|d|=|A\cup B|$. 

 Then an equivalent statement to \eqref{aux} in terms of the space of
$\mathbb{F}_4$-flows $\mathcal{C}_4\cong\mathcal{C}_2\times\mathcal{C}_2$ is that $P(G;4)\neq 0$ if and only if there
exists $z\in(\mathbb{F}_4^\times)^E$ such that
$$\sum_{d\in\mathcal{C}_4,\, d\leq z}(-1)^{|d|}\neq
0.$$ 
This is the assertion of Corollary~\ref{cond Qtensions}.\end{pf}

Theorem~\ref{tripart} is related to the criterion for $P(G;4)\neq 0$ that $G$ be
covered by two bipartite subgraphs $X\cup Y, Y\cup Z$. Given the
latter are bipartite, if
$A\cup C\subseteq X\cup Z$ is eulerian and $C\cup B\subseteq Z\cup Y$
is eulerian then $|A\cup C|$ and
$|C\cup B|$ are both even so that $|A|+|B|$ is also even.
However, a bias in $|A|\!+\!|B|\,(\mbox{\rm mod}\, 2)$ does not imply
a bias in $|A|\!+\!|B|\!+\!|C|\,(\mbox{\rm mod}\, 2)$.

\section{Cubic graphs and triangulations} \label{cubic graphs}

In this section we use MacWilliams duality theorem~\eqref{McW cwe} for
complete weight enumerators to derive a correlation criterion for the existence of 
 a proper vertex
$4$-colouring of a triangulation. 

\begin{thm}\label{4flows tensions} Let $\omega=e^{2\pi i/3}$ and
    let $\psi:\mathbb{F}_4\rightarrow\{0,1,\omega,\omega^2\}$ be a non-trivial
  Dirichlet character (multiplicative, and $\psi(0)=0$). For a graph $G$ with space of
  $\mathbb{F}_4$-flows $\mathcal{C}_4$ and space of $\mathbb{F}_4$-tensions $\mathcal{C}_4^\perp$,
\be\label{first psi}\sum_{z\in\mathcal{C}_4}\prod_{e\in E}\psi(z_e)=2^{n(G)-r(G)}\sum_{z\in\mathcal{C}_4^\perp}\prod_{e\in E}\psi(z_e).\ee
In other words
\be\label{second omega}\mathop{\sum_{\mbox{\rm \tiny eulerian}\, A,B\subseteq E}}_{A\cup
  B=E}\omega^{|A|-|B|}=2^{n(G)-r(G)}\mathop{\sum_{\mbox{\rm \tiny cutsets}\, A,B\subseteq E}}_{A\cup
  B=E}\omega^{|A|-|B|}.\ee

In particular, if $G=(V,E)$ is a cubic graph
 then
\be\label{third cubic}2^{r(G)-n(G)}F(G;4)=\#\{z\in\mathcal{C}_4^\perp:\,\prod_{e\in E}
z_e=1\}-\frac{1}{2}\#\{z\in\mathcal{C}_4^\perp:\,\prod_{e\in E}
z_e\in\{\omega,\omega^2\}\}.\ee
\end{thm}
 
\begin{pf}  We begin by noting that, since
$z\in\mathcal{C}_4$ if and only if $\ol{z}\in\mathcal{C}_4$ and $\psi(\ol{z_e})=\ol{\psi}(z_e)$, 
the equations~\eqref{first psi} and
\eqref{second omega} are between real numbers (in fact rational integers). Also, since $\prod_{e\in E}\psi(z_e)\neq 0$ if
and only if $z$ is nowhere-zero, i.e., $z_e\neq 0$ for all $e\in E$,
the range of the summations in equation~\eqref{first psi} is restricted to nowhere-zero
$\mathbb{F}_4$-flows on the left and nowhere-zero
$\mathbb{F}_4$-tensions on the right.

Identify $\mathbb{F}_4$ with its image
  $\{0,1,\omega,\omega^2\}$ under
  $\psi:\mathbb{F}_4\rightarrow\mathbb{C}$, i.e., $\psi$ is defined by
  $\psi=1_1+\omega 1_\omega+\omega^21_{\omega^2}$. It is easily calculated
  that $\wh{\psi}=2\ol{\psi}$.
By the MacWilliams duality formula~\eqref{McW cwe},
\begin{align*}\sum_{z\in\mathcal{C}_4}\prod_{e\in E}\psi(z_e)={\rm
  cwe}(\mathcal{C}_4;\psi) & =\frac{1}{|\mathcal{C}_4^\perp|}{\rm
  cwe}(\mathcal{C}_4^\perp;\wh{\psi})\\
 & =4^{-r(G)}2^{|E|}{\rm
    cwe}(\mathcal{C}_4^\perp,\ol{\psi}).\end{align*}
Since the
  sums in this equation are real, the function $\ol{\psi}$ can be
  replaced by its conjugate $\psi$, and this establishes equation~\eqref{first psi} of the theorem.

The second
statement~\eqref{second omega} is a straight translation of
\eqref{first psi} into
different language. Using the isomorphism  of
additive groups $\mathbb{F}_4^E\cong\mathbb{F}_2^E\times\mathbb{F}_2^E$,
an element $z\in\mathbb{F}_4^E$ may be written $z=(x,y)$
with $x,y\in\mathbb{F}_2^E$ indicator vectors for subsets
$A,B\subseteq E$ respectively. The property that $z$ is nowhere-zero
translates to the property that $A\cup B=E$ and the condition
$z\in\mathcal{C}_4$ translates to the condition that $A$ and $B$ are
both eulerian. Similarly, the condition
$z\in\mathcal{C}_4^\perp$ translates to the condition that $A$ and $B$
are cutsets.
 
For the final statement~\eqref{third cubic} we use the property that a
cubic graph has a cutset double cover comprising the three-edge stars
at each vertex. For vertex $v\in V$, the three edges $\{e,f,g\}$ incident with $v$
form a star, and each edge occurs exactly twice amongst the $|V|$
stars, since each edge is adjacent to two
distinct vertices.  It follows that $\psi(z_ez_fz_g)\in\{0,1\}$ for
each star $\{e,f,g\}$ and $z\in\mathcal{C}_4$, due to the
 fact that if a sum of three non-zero elements of $\mathbb{F}_4$ is
 equal to $0$ then their product is $1$. 
 Thus we see that
$$\sum_{z\in\mathcal{C}_4}\;\prod_{\mbox{\rm \tiny
    stars}\,\{e,f,g\}}\psi(z_ez_fz_g)=F(G;4).$$
On the other hand, by the double cover property of the collection of stars,
$$\sum_{z\in\mathcal{C}_4}\;\prod_{\mbox{\rm \tiny stars}\,\{e,f,g\}}\psi(z_ez_fz_g)=\sum_{z\in\mathcal{C}_4}\prod_{e\in E}\psi(z_e)^2,$$
and $\psi(z_e)^2=\ol{\psi}(z_e).$
Using equation~\eqref{first psi}, in which $\psi$ is interchangeable
with its conjugate $\ol{\psi}$, this establishes that
$$F(G;4)=2^{n(G)-r(G)}\sum_{z\in\mathcal{C}_4^\perp}\prod_{e\in
  E}\psi(z_e),$$
and equation~\eqref{third cubic} is just another way of writing this.
\end{pf}

We finish this section by interpreting identity~\eqref{third cubic} in
its dual form in terms of the bias of events in a uniform probability
space. The dual notion of a cutset double cover is a cycle double
cover. If $G$ is a plane cubic graph then its planar dual $G^*$ is a plane
triangulation and just as $G$ has a cutset
double cover by three-edge stars (at vertices) so $G^*$ has a
cycle double cover by triangles (faces). 


\begin{thm} \label{triang} Suppose $G$ is a graph that has a cycle double cover by
  triangles and suppose that $\Gamma$ is the event that $A,B\subseteq E$ are eulerian and
$A\cup B=E$. Then, choosing $A,B\subseteq E$
uniformly at random,
 $${\rm Bias}(|A|\equiv |B|\!+\!1\,(\bmod 3)\,\mid\,\Gamma)={\rm
  Bias}(|A|\equiv|B|\!+\!2\,(\bmod 3)\,\mid\,\Gamma)$$
and
 $$\frac{{\rm Bias}(|A|\equiv |B|\,(\bmod 3)\,\mid\,\Gamma)-{\rm
    Bias}(|A|\equiv |B|\!+\!1\,(\bmod 3)\,\mid\,\Gamma)}{{\rm Bias}(|A|\equiv |B|\,(\bmod 3)\,)-{\rm
  Bias}(|A|\equiv|B|\!+\!1\,(\bmod 3)\,)}=\frac{2^{3|E|-2|V|}P(G;4)}{F(G;4)}.$$ 
\end{thm}
\begin{pf}  
The left-hand sum in equation~\eqref{second omega} of Theorem~\ref{4flows tensions} has
the following interpretation:
\be\label{as probability}2^{-2|E|}\mathop{\sum_{\mbox{\rm \tiny eulerian}\, A,B\subseteq E}}_{A\cup
  B=E}\omega^{|A|-|B|}\ee
\begin{align*} & =\mathbb{P}(|A|\equiv|B|\,(\bmod\,3)\,\cap\,\Gamma)+\omega\mathbb{P}(|A|\equiv|B|\!+\!1\,(\bmod\, 3)\,\cap\,\Gamma)+\omega^2\mathbb{P}(|A|\equiv|B|\!+\!2\,(\bmod\, 3)\,\cap\,\Gamma)\\
 & = \mathbb{P}(|A|\equiv|B|\,(\bmod\,
  3)\,\cap\,\Gamma)+\omega^2\mathbb{P}(|A|\equiv|B|\!+\!1\,(\bmod\,
  3)\,\cap\,\Gamma)+\omega\mathbb{P}(|A|\equiv|B|\!+\!2\,(\bmod\, 3)\,\cap\,\Gamma),\end{align*}
the latter equality since, as remarked in the proof of
Theorem~\ref{4flows tensions}, the sum we started with is
  real.
Hence 
\be\label{1 equals 2}\mathbb{P}(|A|\equiv|B|\!+\!1\,(\bmod\,
3)\,\cap\,\Gamma)=\mathbb{P}(|A|\equiv|B|\!+\!2\,(\bmod\, 3)\,\cap\,\Gamma).\ee
By definition of $\Gamma$ and since $G$ has a cycle double cover by triangles,
$\mathbb{P}(\Gamma)=2^{-2|E|}F(G;4)\neq 0$.
Dividing equation~\eqref{1 equals 2} by $\mathbb{P}(\Gamma)$ yields the first statement of the theorem. 

Equation~\eqref{1 equals 2} together with the identity developed in
\eqref{as probability} has the consequence that, in the notation of
Theorem~\ref{4flows tensions},
\be\label{lhs as prob}2^{-2|E|}\sum_{z\in\mathcal{C}_4}\prod_{e\in
  E}\psi(z_e)=\mathbb{P}(|A|\equiv|B|\,(\bmod\,
3)\,\cap\,\Gamma)-\mathbb{P}(|A|\equiv|B|\!+\!1\,(\bmod\, 3)\,\cap\,\Gamma).\ee

 By equation~\eqref{first
  psi} of Theorem~\ref{4flows tensions}, in which $\prod_{e\in
  E}\psi(z_e)\in\{0,1\}$ for $z\in\mathcal{C}_4^\perp$ since $G$ has a
cycle double cover by triangles, and equation~\eqref{lhs as prob},
\begin{align*}\mathbb{P}(|A|\equiv|B|\,(\bmod 3)\,\cap\,\Gamma)-\mathbb{P}(|A|\equiv|B|\!+\!1\,(\bmod
3)\,\cap\,\Gamma) & =
2^{-2|E|}2^{n(G)-r(G)}\sum_{z\in\mathcal{C}_4^\perp}\prod_{e\in E}\psi(z_e)\\
 & = 2^{-2|E|}2^{n(G)-r(G)}4^{-k(G)}P(G;4).
\end{align*}
Dividing this last equation through by $\mathbb{P}(\Gamma)=2^{-2|E|}F(G;4)$ gives
$$\mathbb{P}(|A|\equiv|B|\,(\bmod
3)\,\mid\,\Gamma)-\mathbb{P}(|A|\equiv|B|\!+\!1\,(\bmod
3)\,\mid\,\Gamma)=2^{|E|-2|V|}P(G;4)/F(G;4),$$
i.e.,
$${\rm Bias}(|A|\equiv |B|\,(\bmod 3)\,\mid\,\Gamma)-{\rm
    Bias}(|A|\equiv |B|\!+\!1\,(\bmod 3)\,\mid\,\Gamma)=2^{|E|-2|V|+1}P(G;4)/F(G;4)$$

By Lemma~\ref{bias S} with $q=3$ and
  $\Sigma=\{A,B\subseteq E: |A|-|B|\in S\}$, taking the difference
  between the cases $S=\{0\}$ and $S=\{1\}$ we obtain
$${\rm Bias}(|A|\!-\!|B|\equiv 0\,(\bmod\, 3)\,)-{\rm
  Bias}(|A|\!-\!|B|\equiv 1\,(\bmod\, 3)\,)$$
\begin{align*} & =3^{-1}2^{1-|E|}\big[(1-\mbox{$\frac{1}{2}$})^{|E|}(1-e^{2\pi i/3})+(1-\mbox{$\frac{1}{2}$})^{|E|}(1-e^{4\pi i/3})\big]\\
 & =2^{1-2|E|}\end{align*}
 The second statement of the theorem now results.
\end{pf}

In particular, a graph $G$ with a cycle double cover by triangles has $P(G;4)\neq 0$ if and only if
$\mathbb{P}(|A|\equiv|B|\,(\mbox{\rm mod}\, 3)\,\mid\,
\Gamma)>\frac{1}{3}$, i.e., the event that $A,B$ form an eulerian cover
of $G$ is positively correlated with $|A|\equiv|B|\,(\bmod\, 3)$.

\section{Eulerian subdigraphs of a $4$-regular graph}\label{4 regular}


In this final section we use MacWilliams duality \eqref{McW cwe} for
compete weight enumerators and Lemma~\ref{cwe sum} to derive some
further evaluations of the Tutte polynomial on $H_3$ similar in form
to Theorem~\ref{triang}. 

Take $q=3$ and $G$ a $4$-regular graph (such as the line graph of
a plane cubic graph), for which  the space of $\mathbb{F}_3$-flows has a natural
identification with the set of eulerian partial orientations of $G$. In a partial
orientation of a graph some edges may not be directed; in an {\em 
eulerian} partial orientation each vertex
has the same number of incoming and outgoing directed edges.
 
A reference orientation $\gamma$ of $G$ is fixed. A partial
orientation $\alpha$ is defined corresponding to a
vector $x\in\mathbb{F}_3^E$ by making $\alpha$ direct an
edge $e$ the same way as $\gamma$ if $x_e=+1$, making $\alpha$ reverse the
direction of $\gamma$ if $x_e=-1$, and leaving $e$ undirected if
$x_e=0$.
Given that $G$ is $4$-regular, if $x$ is a $\mathbb{F}_3$-flow 
then it defines an eulerian partial orientation $\alpha$.    
If further $x$ is nowhere-zero then $\alpha$ is an eulerian orientation of
$G$.
For orientations $\alpha,\beta$, define $\alpha+\beta$ to be
the partial orientation
whose directed edges are those sharing the same direction in
$\alpha$ and $\beta$. If $x,y\in\mathbb{F}_3^E$ define $\alpha,\beta$
relative to the base orientation $\gamma$ of $G$ then $\alpha+\beta$
is the partial orientation defined by $-(x+y)$.

Suppose an orientation $\alpha$ is chosen uniformly at random. Let $\Sigma$ be the event that $|\alpha+\gamma|\equiv 0\,(\mbox{\rm
  mod}\, 2)$. 
Clearly ${\rm Bias}(\Sigma)=0.$
Let $\Gamma$ be the event that $\alpha$ is an eulerian orientation.
Then $\mathbb{P}(\Gamma)=2^{-|E|}F(G;3)$, since for the $4$-regular
graph $G$ the number of eulerian orientations is the number of
nowhere-zero $3$-flows. 

Finding ${\rm Bias}(\Sigma\,\mid\,\Gamma)$ is a bit more difficult and
in order to state a partial result on this we need some definitions.
The {\em line graph} $L(H)$ of a graph $H$ has vertices the edges of $H$ and
adjacent vertices $e,f$ when $e$ and $f$ are incident in $H$. If $H$
is embedded in an orientable surface, the {\em medial graph} $M(H)$ of
$H$ is the graph obtained by placing vertices at the edges
of $H$ and joining vertices $e,f$ of $M(H)$ by an edge if they lie on
incident edges $e,f$ of $H$ and it is possible to draw a line joining
$e$ and $f$ without crossing any edges of $H$. (If edges $e,f$ are incident with
a vertex of degree $2$ then they are joined by two edges, neither of
which can be continuously transformed to the other without crossing
an edge of $H$.)
The medial graph $M(H)$ is $4$-regular. 

Suppose now that $H$ is an
orientably embedded cubic graph. Then $M(H)$
is an embedding of $L(H)$ in the same orientable surface as $H$. A
vertex of $H$ lies in the interior of a
triangle of edges in $M(H)$, which we shall call a {\em black triangle} of
$M(H)$ (on account of the standard white-black face colouring of the medial graph). When the edges of the black triangles of $M(H)$ are directed in a
clockwise sense on the surface on which $M(H)$ is embedded, the
resulting orientation of $M(H)$ is eulerian. (The clockwise direction
traced by the edges of
a black triangle of $M(H)$ corresponds to a clockwise orientation of the
three edges at a vertex of $H$, called a vertex rotation in the
embedding of $H$.)

\begin{thm}\label{even odd eul suborns} Let $G$ be the medial graph of a plane cubic
  graph and let $\gamma$ be the orientation directing edges of $G$
  clockwise around black triangles. Then, choosing an orientation
  $\alpha$ of $G$ uniformly at random, the event $\Sigma$ that
  $\alpha$ agrees with $\gamma$ on an even number of edges and the
  event $\Gamma$ that $\alpha$ is eulerian have correlation given by 
$${\rm Bias}(\Sigma\,\mid\,\Gamma)=\frac{P(G;3)}{F(G;3)}.$$
\end{thm}
\begin{pf} Let $\mathcal{C}_3$ be the space of $\mathbb{F}_3$-flows of
  $G$. Then
\begin{align*}\mathbb{P}(\Gamma){\rm Bias}(\Sigma\,\mid\,\Gamma) & =\mathbb{P}(\Sigma\cap\Gamma)-\mathbb{P}(\ol{\Sigma}\cap\Gamma)\\
  & =2^{-|E|}{\rm cwe}(\mathcal{C}_3;1_1-1_{-1})\\
 & = 2^{-|E|}3^{-r(G)}{\rm
   cwe}(\mathcal{C}_3^\perp;(-3)^{\frac{1}{2}}(1_{-1}-1_{1})\\
& = (-1)^{|V|}2^{-|E|}3^{k(G)}{\rm cwe}(\mathcal{C}_3^\perp;1_1-1_{-1}),\end{align*}
using $|E|=2|V|$ for $4$-regular graph $G$. A result of Penrose\footnote{Penrose quoted this theorem (in a different
  formulation) in \cite{P71}, in which he mentioned that his proof was too lengthy for inclusion. An
  elegant short proof has been given by Kaufmann~\cite{K90}. See also
  \cite[Theorem 3.1]{EG96} for a
  generalisation and for further citations  -
  for example Scheim~\cite{S74} found the result
  independently and was the first to
  publish a proof. 

It remains an open problem~\cite{EG96} to characterise those edge $3$-colourable cubic graphs $H$ for
which the line graph $L(H)$ with a fixed orientation of its edges has the property that the
number of nowhere-zero $\mathbb{F}_3$-tensions of $L(H)$ with an even number of
edges with value $-1$ differs from those with an odd number of edges
with value $-1$. The line graph $L(K_{3,3})$ does not have this
property. The theorem of Penrose and Scheim is that when $H$ is a planar cubic
graph nowhere-zero $\mathbb{F}_3$-tensions of $L(H)$ either all have an
 even number of edges with value $-1$ or all an odd number of such
 edges.} says that
for any nowhere-zero $\mathbb{F}_3$-tension $x$ of $G$ (corresponding to
an edge $3$-colouring of the plane cubic graph $H$) we have
$$\prod_{e\in E}(1_1-1_{-1})(x_e)=(-1)^{\#\{e\in E:x_e=-1\}}=(-1)^{|V|}$$ 
when $G$ has its fixed orientation $\gamma$ clockwise around black triangles (or any
other orientation $\beta$ with $|\beta+\gamma|$ even). With
$\mathbb{P}(\Gamma)=2^{-|E|}F(G;3)\neq 0$ the theorem is proved.
\end{pf}

We move on now to choosing pairs of orientations of a $4$-regular graph
and formulate an analogue of Proposition~\ref{squares with q=2} for
eulerian orientations rather than the eulerian subgraphs of that proposition.

Suppose $\alpha,\beta$ are orientations of $G$ chosen uniformly at random.
Let $\Sigma$ be the event that $|\alpha+\beta|$ is even. 
Note that if $\beta$ is another orientation then
  $|\alpha+\gamma|$ and $|\beta+\gamma|$ have the same parity if and
  only if $|\alpha+\beta|$ is even. We have ${\rm Bias}(\Sigma)=0$ as before.
Let $\Gamma$ be the event that $\alpha+\beta$ is an eulerian partial
orientation of $G$, i.e., the consistently directed edges of $\alpha$
and $\beta$ form an eulerian partial orientation of $G$. 

\begin{lem} \label{PGamma} {\rm (Cf. \cite[Corollary 5]{AJG05})}
  The event $\Gamma$ that orientations $\alpha,\beta$ of a
  $4$-regular graph agree on an eulerian partial orientation has
  probability
$$\mathbb{P}(\Gamma)=4^{-|E|}T(G;2,4).$$
\end{lem}
\begin{pf} Let $\mathcal{C}_3$ be the space of $\mathbb{F}_3$-flows of $G$.
If orientations $\alpha,\beta$ are defined relative to the base
orientation $\gamma$ of $G$ by the nowhere-zero vectors
$x, -y\in\mathbb{F}_3^E$, then the event that $\alpha+\beta$ is eulerian
coincides with the event that $x-y\in\mathcal{C}_3$, i.e., $x,y$ belong
to the same coset of $\mathcal{C}_3$.  Using Lemma~\ref{cwe sum} 
\begin{align*}\mathbb{P}(\Gamma) &
  =4^{-|E|}\sum_{\mathcal{C}_3+z\in\mathbb{F}_3^E/\mathcal{C}_3}{\rm cwe}(\mathcal{C}_3+z;1_{1,-1})^2\\
 & = 4^{-|E|}3^{-r(G)}{\rm cwe}(\mathcal{C}_3^\perp;(21_0-1_{1,-1})^2)\\
 & = 4^{-|E|}3^{-r(G)}{\rm hwe}(\mathcal{C}_3^\perp;4)\\
 & =4^{-|E|}T(G;2,4).\end{align*}
\end{pf} 

This lemma leads us to our promised theorem.

\begin{thm}\label{eul orn correlation}
 Let $G$ be a $4$-regular graph and $\alpha,\beta$ two orientations of
  $G$ chosen uniformly at random. Then the event $\Sigma$ that
  $\alpha$ agrees with $\beta$ on an even number of edges and the
  event $\Gamma$ that $\alpha$ agrees with $\beta$ in an eulerian partial
  orientation of $G$ have correlation given by 
$${\rm Bias}(\Sigma\,\mid\,\Gamma)=\frac{3^{|E|-|V|}P(G;3)}{T(G;2,4)}.$$
\end{thm} 
\begin{pf} Using Lemma~\ref{cwe sum} 
\begin{align*}\mathbb{P}(\Gamma){\rm Bias}(\Sigma\,\mid\,\Gamma) & =
  \mathbb{P}(\Sigma\cap\Gamma)-\mathbb{P}(\ol{\Sigma}\cap\Gamma)\\
 & = 4^{-|E|}\sum_{\mathcal{C}_3+z\in\mathbb{F}_3^E/\mathcal{C}_3}|{\rm
   cwe}(\mathcal{C}_3+z;1_1-1_{-1})|^2\\
 & = 4^{-|E|}3^{-r(G)}{\rm
   cwe}(\mathcal{C}_3^\perp;|(-3)^{\frac{1}{2}}(1_{-1}-1_1)|^2)\\
 & =4^{-|E|}3^{-r(G)}{\rm
   cwe}(\mathcal{C}_3^\perp;31_{1,-1})=4^{-|E|}3^{n(G)}(-1)^{r(G)}T(G;-2,0)\\
 & = 4^{-|E|}3^{|E|-|V|}P(G;3).\end{align*}
Lemma~\ref{PGamma} now gives the result.
\end{pf}


 
{\small 
\bibliography{4flow.bib} 

\begin{thebibliography}{10}

\bibitem{A99}
N.~Alon.
\newblock {Combinatorial Nullstellensatz}.
\newblock {\em Combinatorics, Probability and Computing}, 8(1 \& 2):7--29,
  1999.

\bibitem{AT92}
N.~Alon and M.~Tarsi.
\newblock Colorings and orientations of graphs.
\newblock {\em Combinatorica}, 12:125--134, 1992.

\bibitem{AT97}
N.~Alon and M.~Tarsi.
\newblock A note on graph colorings and graph polynomials.
\newblock {\em Journal of Combinatorial Theory Series B}, 70:197--201, 1997.

\bibitem{B98}
B.~Bollob\'{a}s.
\newblock {\em Modern Graph Theory}.
\newblock Springer, New York, 1998.

\bibitem{EG96}
M.~N. Ellingham and L.~Goddyn.
\newblock List edge colourings of some 1-factorable multigraphs.
\newblock {\em Combinatorica}, 16:343--352, 1996.

\bibitem{E98}
R.~Eriksson.
\newblock A property of coset weight distributions.
\newblock In {\em Proceedings of the sixth international workshop on algebraic
  and combinatorial coding theory}, pages 112--115, Pskov, Russia, 1998.

\bibitem{AJG05}
A.~Goodall.
\newblock {Some new evaluations of the Tutte polynomial}.
\newblock {\em Journal of Combinatorial Theory Series B}, 96:207--224, 2006.

\bibitem{CCC07}
A.~Goodall.
\newblock Fourier analysis on finite abelian groups: some graphical
  applications.
\newblock In G.~Grimmett and C.~McDiarmid, editors, {\em {Combinatorics,
  Complexity, Chance: A Tribute to Dominic Welsh}}, number~34 in Oxford Lecture
  Series in Mathematics and its Applications, chapter~7, pages 103--129. Oxford
  University Press, Oxford, 2007.

\bibitem{K90}
L.~Kaufmann.
\newblock Map coloring and the vector cross product.
\newblock {\em Journal of Combinatorial Theory Series B}, 48:145--154, 1990.

\bibitem{Kung07}
J.~Kung.
\newblock Coboundary, flows and {T}utte polynomials of matrices.
\newblock {\em Annals of Combinatorics}.
\newblock To appear.

\bibitem{M04}
Y.~Matiyasevich.
\newblock Some probabilistic restatements of the four color conjecture.
\newblock {\em Journal of Graph Theory}, 46:167--179, 2004.

\bibitem{Mi05}
O.~Milenkovic.
\newblock Support weight enumerators and coset weight distributions of isodual
  codes.
\newblock {\em Designs, Codes and Cryptography}, 35(1):81--109, 2005.

\bibitem{Onn04}
S.~Onn.
\newblock Nowhere-zero flow polynomials.
\newblock {\em Journal of Combinatorial Theory Series $A$}, 108:205--215, 2004.

\bibitem{P71}
R.~Penrose.
\newblock Applications of negative dimensional tensors.
\newblock In D.~Welsh, editor, {\em Combinatorial Mathematics and its
  Applications}, pages 221--244, New York, London, 1971. Academic Press.

\bibitem{RR78}
R.~C. Read and P.~Rosenstiehl.
\newblock On the principal edge tripartition of a graph.
\newblock {\em Annals of Discrete Mathematics}, 3:195--226, 1978.

\bibitem{S74}
D.~E. Scheim.
\newblock The number of edge $3$-colorings of a planar cubic graph as a
  permanent.
\newblock {\em Discrete Mathematics}, 8:377--382, 1974.

\bibitem{St97}
R.~Stanley.
\newblock {\em Enumerative Combinatorics}, volume~1.
\newblock Cambridge University Press, 1997.

\bibitem{Tarsi}
M.~Tarsi.
\newblock The graph polynomial and the number of proper vertex colorings.
\newblock {\em Annales de l'Institut Fourier}, 49(3):1089--1093, 1999.

\bibitem{Terras99}
A.~Terras.
\newblock {\em Fourier analysis on finite groups and applications}.
\newblock Cambridge University Press, 1999.

\bibitem{vdW41}
B.~van~der Waerden.
\newblock {Die lange Reichweite der regelm\"{a}ssigen Atomanordnung in
  Mischkristallen}.
\newblock {\em Z. Physik}, 118:473, 1941.

\bibitem{DW93}
D.~Welsh.
\newblock {\em Complexity: {K}nots, {C}olourings and {C}ounting}.
\newblock Cambridge University Press, 1993.

\end{thebibliography}
\bibliographystyle{abbrv} 
}

\end{document}